\documentstyle{amsppt}
\magnification=1200
\hoffset=-0.5pc
\vsize=57.2truepc
\hsize=38truepc
\nologo
\spaceskip=.5em plus.25em minus.20em
\define\fra{\frak}
\define\Bobb{\Bbb}
\define\rh{r}
\define\wide{}
 \define\atiyboo{1}
 \define\atibottw{2}
\define\bodenhu{3}
\define\borelone{4}
\define\broediec{5}
\define\brownboo{6}
\define\guhujewe{7}
\define\helgaboo{8}
   \define\modus{9}
\define\modustwo{10}
    \define\srni{11}
\define\huebjeff{12}
\define\jeffrone{13}
\define\jeffrtwo{14}
\define\karshone{15}
\define\kirilboo{16}
\define\mehtsesh{17}
\define\narasesh{18}
\define\rasheone{19}
\define\sjamlerm{20}
\define\trottone{21}
 \define\weiltwo{22}
\define\weinsthi{23}
\noindent
dg-ga/9710033

\topmatter
\title 
On the variation of the Poisson structures\\
of certain moduli spaces
\endtitle
\author 
Johannes Huebschmann
\endauthor
\affil 
Universit\'e des Sciences et Technologies
\\
UFR de Math\'ematiques
\\
F-59 655 Villeneuve d'Ascq C\'edex
\\
France
\\
Johannes.Huebschmann\@univ-lille1.fr
\endaffil
\abstract{Given a Lie group G whose Lie algebra is endowed with a nondegenerate
invariant symmetric bilinear form, we construct a Poisson algebra of continuous
functions on a certain open subspace $\Cal R(\pi,G)$ of the space of 
representations in G of the fundamental group $\pi$ of a compact connected 
orientable topological surface with finitely many boundary circles; when G is 
compact and connected, $\Cal R(\pi,G)$ may be taken dense in the space of all 
representations. The space $\Cal R(\pi,G)$ contains spaces of representations 
where the values of those generators of the fundamental group which correspond 
to the boundary circles are constrained to lie in fixed conjugacy classes and, 
on these representation spaces, the Poisson algebra restricts to stratified 
symplectic Poisson algebras constructed elsewhere earlier. Hence the Poisson 
algebra on $\Cal R(\pi,G)$ gives a description of the variation of the 
stratified symplectic Poisson structures on the smaller representation spaces 
as the chosen conjugacy classes move.} 
\endabstract
\endkeywords
\endsubjclass
\endtopmatter
\document
\rightheadtext{Variation of Poisson structure}

\beginsection 1. Introduction

Given a Lie group $G$ 
whose Lie algebra is endowed with a nondegenerate
invariant symmetric  bilinear form, 
certain subspaces $\roman{Rep}(\pi,G)_{\bold C}$
(the subscript will be explained below)
of the space $\roman{Rep}(\pi,G)$
of representations in $G$ of the fundamental group 
$\pi$ of an orientable topological
surface with finitely many boundary circles
are known
to inherit
symplectic or more generally Poisson structures;
see e.~g. \cite \guhujewe\ and the literature there.
These subspaces
are obtained
when the values of those generators 
of $\pi$ which correspond to the boundary circles 
are constrained to lie in fixed conjugacy classes of $G$;
apart from trivial cases, these
subspaces have stricly positive codimension 
in $\roman{Rep}(\pi,G)$ and hence are genuinely smaller than
$\roman{Rep}(\pi,G)$.
But what can be said about the {\it variation\/} of the
{\it Poisson structures\/} when the 
{\it conjugacy classes\/} are allowed to {\it move\/}?
It is the objective of the present paper to offer an answer
to this question, phrased in the world of
{\it Poisson geometry\/}.
Our main result will give a Poisson structure
on a suitable subspace $\Cal R(\pi,G)$
of $\roman{Rep}(\pi,G)$
which is a union 
of distinct spaces of the kind $\roman{Rep}(\pi,G)_{\bold C}$;
this Poisson structure will induce the Poisson structures
on the spaces $\roman{Rep}(\pi,G)_{\bold C}$.
When $G$ is compact and connected,
$\Cal R(\pi,G)$ may in fact be taken to be dense in $\roman{Rep}(\pi,G)$.
We now explain this in some more detail.
\smallskip
Let $\Sigma$ be a compact connected orientable topological surface
of genus $\ell$
with boundary $\partial \Sigma$ consisting of $n$ circles
$S_1,\dots,S_n$.
To avoid trivial cases or inconsistencies,
when the genus is zero,
we suppose that $n \geq 3$.
Consider the usual 
presentation
$$
\Cal P = \langle x_1,y_1,  \dots,x_\ell,y_\ell,z_1,\dots,z_n; r
\rangle,
\quad
r = \Pi_{j=1}^{\ell} [x_j,y_j] z_1 \dots z_n,
\tag1.1
$$
of the fundamental group $\pi = \pi_1(\Sigma)$.
Let $\cdot$ be a nondegenerate invariant symmetric bilinear form
on the Lie algebra $\fra g$ of $G$.
Extending the approach in \cite\guhujewe,
from these data, we shall  construct a Poisson structure
on a certain open  subspace 
$\Cal R(\pi,G)$
of 
the space $\roman {Rep}(\pi,G)$ 
of all representations of $\pi$
in $G$, $\roman {Rep}(\pi,G)$  being
the orbit space for the action of $G$ on
the space $\roman{Hom}(\pi,G)$ of homomorphisms
by conjugation.
\smallskip
Given
an $n$-tuple 
$\bold C = (C_1,\dots,C_n)$
of conjugacy classes in $G$,
denote by
$\roman{Hom}(\pi,G)_{\bold C}$
the space of homomorphisms $\chi$ from $\pi$ to $G$
for which the value $\chi (z_k)$
of each generator $z_k$
lies in $C_k$, for $1 \leq k \leq n$,
and denote by
$\roman{Rep}(\pi,G)_{\bold C}$
the corresponding space of representations,
that is, the orbit space for the action of $G$ on
$\roman{Hom}(\pi,G)_{\bold C}$
by conjugation.
From these
data, 
in \cite\guhujewe,
we constructed
the structure of a stratified symplectic space
on $\roman{Rep}(\pi,G)_{\bold C}$, in particular,
a Poisson algebra
$\left(C^{\infty}(\roman{Rep}(\pi,G)_{\bold C}),\{\cdot,\cdot\}\right)$
of continuous functions
on $\roman{Rep}(\pi,G)_{\bold C}$; the latter endows a certain top stratum
with the structure of a symplectic manifold.
Henceforth we shall refer to a Poisson structure
which
is part of the  structure of a stratified symplectic space
as a {\it stratified symplectic Poisson structure\/};
thus a stratified symplectic Poisson algebra on a space with a
single stratum is just a smooth symplectic Poisson algebra
in the usual sense.
\smallskip
Let $O$ be
an open connected Ad-invariant neighborhood of
zero in $\fra g$ which, under the exponential map
exp from $\fra g$ to $G$, is mapped {\it diffeomorphically\/}
onto
an open invariant neighborhood $B$ of
the identity of $G$.
When $G$ is compact,  connected and simply connected,
and hence semisimple,
a maximal such subset $O$ is given by
the largest {\it connected\/} neighborhood of the origin
of $\fra g$ where the exponential map is regular;
it consists of those $X$ in $\fra g$
which have the property that the endomorphism $\roman{ad}(X)$
of $\fra g$ has only eigenvalues $\lambda = 2 \pi i \nu$
with $|\nu| <1$;
see Lemma 3.1 and Corollary 3.3 below.
The image $B$ in $G$ of $O$ then consists precisely of the regular values
of the exponential map,
and this choice of a maximal 
connected neighborhood $O$ of
zero in $\fra g$
is in fact canonical;
furthermore, $B$ then contains all regular elements of $G$
in the sense of Lie groups 
(they correspond to the stronger condition
$0 < |\nu| <1$)
whence $B$ is certainly dense in $G$.
For example, for $G=\roman {SU}(2)$, we then have
$B = G \setminus \{-\roman{Id}\}$.
What happens for 
$G=\roman {SU}(n)$ when $n \geq 3$ will be briefly explained in (3.5) below.
On the other hand,
when $G$ is still compact, connected and semisimple
but not simply connected, things get more complicated;
for example,
for $G=\roman {SO}(3)$, 
an appropriate subset $B$ is the ball in
$\roman {SO}(3)$ which arises when the conjugacy class of
$\left[\matrix -1& 0 & 0\\
                0&-1 & 0\\
                0& 0 & 1
       \endmatrix
 \right]$
is removed from 
$\roman {SO}(3)$;
this conjugacy class is a copy of the real projective plane,
and its points are still regular values for the exponential map.
Returning to the general case,
let $\Cal H(\pi,G)$ be
the space of homomorphisms $\chi$ from $\pi$ to $G$
for which the value $\chi (z_k)$
of each generator $z_k$ lies in $B$
and let
$\Cal R(\pi,G)$ be
the corresponding space of $G$-orbits;
we do not indicate the dependence on $B$.
We shall say that a conjugacy class $C$ is $B$-{\it regular\/}
if it is contained in $B$.
It is manifest that
$\Cal R(\pi,G)$
is the union in
$\roman{Rep}(\pi,G)$
of the 
spaces $\roman{Rep}(\pi,G)_{\bold C}$
for those 
$n$-tuples 
$\bold C = (C_1,\dots,C_n)$
of conjugacy classes which are $B$-{\it regular\/} in the sense
that
every $C_k$ is $B$-regular, and 
$\Cal R(\pi,G)$
is plainly open 
in $\roman{Rep}(\pi,G)$.

\proclaim{Theorem}
The data determine the structure of a Poisson algebra of 
continuous functions
on 
$\Cal R(\pi,G)$
which, on each 
moduli space $\roman{Rep}(\pi,G)_{\bold C}$ (that
lies in $\Cal R(\pi,G)$),
restricts to the 
stratified symplectic 
Poisson algebra
on $\roman{Rep}(\pi,G)_{\bold C}$
determined by the data.
When $G$ is compact and connected, 
the subspace
$\Cal R(\pi,G)$ may be taken
to be dense in $\roman{Rep}(\pi,G)$.
\endproclaim

In particular,
the strata of the 
spaces $\roman{Rep}(\pi,G)_{\bold C}$
in $\Cal R(\pi,G)$,
which carry symplectic structures in view of what has been said earlier,
appear
as symplectic leaves (in a generalized sense)
for the Poisson structure on
$\Cal R(\pi,G)$.
The first statement of the theorem will be proved in (2.18) below,
and the second statement, 
which says that, for $G$
compact and connected, 
the subspace
$\Cal R(\pi,G)$ may be taken to be
dense in $\roman{Rep}(\pi,G)$,
will be justified in Section 3.
\smallskip
We now explain briefly how this theorem is proved.
By {\it purely finite dimensional methods\/},
we shall construct a suitable {\it extended moduli space\/},
that is to say, a smooth symplectic manifold $\Cal M$
and a hamiltonian action of
the group 
$G_0 \times G_1\times \dots \times G_n$ on $\Cal M$,
each $G_j$ being a copy of $G$, $0 \leq j \leq n$,
with momentum mapping $\mu_j \colon \Cal M \to \fra g^*_j$
for the action of the $j$'th copy  $G_j$ of $G$,
in such a way that
$$
\mu_0 \times \mu_1 \times \dots\times \mu_n
\colon
\Cal M 
@>>>
\fra g_0^* \times \fra g_1^*\times \dots \times \fra g_n^*
$$
is a momentum mapping for the
action of
$G_0 \times G_1\times \dots \times G_n$ on $\Cal M$
and that the following holds:
Consider
the reduced space 
$\Cal M_0 = \mu_0^{-1}(0)\big/ G$
for the hamiltonian action
of the zero'th copy  $G_0$ of $G$ on $\Cal M$,
with its structure 
$\left(C^{\infty}(\Cal M_0),\{\cdot,\cdot\}\right)$
of a stratified symplectic space
\cite\sjamlerm.
The action
of $G_0 \times G_1\times \dots \times G_n$ on $\Cal M$
induces a hamiltonian action 
(in the sense of stratified symplectic spaces)
of
$\Gamma =G_1\times \dots \times G_n$
on $\Cal M_0$,
and there is a canonical map
from $\Cal M_0$
to $\Cal R(\pi,G)$
which induces a homeomorphism
between the space of orbits
$\Cal M_0\big/ \Gamma$
and
$\Cal R(\pi,G)$.
The 
subalgebra 
$(C^{\infty}(\Cal M_0))^\Gamma$
of $\Gamma$-invariant
functions in 
$C^{\infty}(\Cal M_0)$
then yields a Poisson algebra 
$\left(C^{\infty}(\Cal R(\pi,G)),\{\cdot,\cdot\}\right)$
of continuous functions
on $\Cal R(\pi,G)$
having the asserted properties.
\smallskip
The philosophy is the same as that in \cite\guhujewe:
The concept of fundamental group is too weak to handle
peripheral structures;
in that paper,
the fundamental group has been replaced by two more general concepts
which enabled us to overcome the difficulties with the peripheral structure:
by that of a {\it group system\/},
to handle the {\it global structure\/}, and by that of a
{\it fundamental groupoid\/},
to handle the infinitesimal structure.
In the present paper, we shall 
study the {\it global\/} structure by means of
the {\it fundamental groupoid\/}. 
This leads to 
the extended moduli space
$\Cal M$ mentioned above;
this space is somewhat \lq\lq larger\rq\rq\ 
than the extended moduli space in \cite\guhujewe\ in the sense
that the extended moduli space
in \cite\guhujewe\ 
which arises from a choice of  conjugacy classes
that are $B$-regular
is obtained by symplectic reduction 
from the extended moduli space
$\Cal M$ to be explored in the present paper,
see (2.14) below;
on the other hand, 
those conjugacy classes  which are not $B$-regular, whatever choice of $B$,
{\it cannot\/} be incorporated
in the approach in terms of this larger
extended moduli space, though.
We shall explain the obstacles in
Remark 2.19 below.
\smallskip
Write
$\Sigma^{\bullet}$ for the corresponding punctured surface,
which contains $\Sigma$ and,
for each boundary component
of $\Sigma$,
has a single puncture
in such a way that,
for $1 \leq k \leq n$,
the $k$'th boundary circle $S_k$ surrounds the $k$'th puncture
and no other puncture.
For $G = \roman U(r)$, the unitary group,
after a choice of complex
structure on
$\Sigma^{\bullet}$
has been made,
a space of the kind
$\roman{Rep}(\pi,G)_{\bold C}$
admits an interpretation as moduli space
of semistable 
rank $r$
parabolic bundles
of parabolic degree zero,
with flags and weights
given by $\bold C$.
In fact,
the search for a sensible moduli space
of rank $r$ holomorphic vector bundles
on
$\Sigma^{\bullet}$
in the
category of complex varieties
led to the
introduction of the additional structure
of flags and weights at the punctures
of $\Sigma^{\bullet}$
\cite\mehtsesh.
Now, whether or not $G$ is the unitary group,
in our representation space picture,
a complex structure of $\Sigma^{\bullet}$
does not come into play and
the complex structure on the
moduli space
is not visible
but, from the results of \cite\guhujewe,
$\roman{Rep}(\pi,G)_{\bold C}$
comes with a structure of a stratified symplectic space.
The main result of the present paper,
the theorem spelled out above,
gives a Poisson structure
on a certain ambient space
$\Cal R(\pi,G)$
which,
among other things,
yields
a {\it description
of the variation of the stratified symplectic structures
across the pieces
$\roman{Rep}(\pi,G)_{\bold C}$
of
$\Cal R(\pi,G)$
as $\bold C$
moves\/}.
For $G$ the unitary group,
the variation of the complex structure
of the parabolic moduli spaces
as the flags and weights
change has been
studied in
\cite\bodenhu,
but {\it not\/} by means of an ambient space.
\smallskip
We keep
the notations of \cite\guhujewe,
and most unexplained concepts may be found in that paper.
A leisurely introduction into Poisson geometry
of certain moduli spaces
(related to those
studied in the present paper
but not exactly the same ones)
may be found in \cite\srni\  where also more references are given.
A  space similar to our extended moduli space
$\Cal M$
has been given in \cite\jeffrone\ and, by infinite dimensional methods
from gauge theory,
a symplectic structure on that space has been constructed.
See Remark 2.20 below for details.
I am much indebted to
A. Weinstein for discussions at various stages;
he was presumably the first to notice the significance of
groupoids for moduli spaces.
Thanks are also due to S. Helgason and R. Steinberg for having
pointed out to me the relevant literature where the simple connectedness
of conjugacy classes of compact, connected, and simply connected
Lie groups has been established.

\beginsection 2. Details and Proofs

The system
$(\pi; \pi_1,\dots, \pi_n)$,
with
$\pi = \pi_1(\Sigma)$ and
$\pi_k = \pi_1(S_k) \cong \bold Z$,
for $1 \leq k \leq n$,
is what is called a {\it group system\/}
\cite{\guhujewe,\ \trottone}.
When the boundary $\partial \Sigma$
of $\Sigma$
is non-empty
the group $\pi$ is free; yet it is convenient to use the
presentation (1.1).
\smallskip
Next we recall
the appropriate fundamental groupoid:
Pick a base point $p_0$ 
not on the boundary
and, moreover,
for each boundary component $S_k$ of $\Sigma$,
pick a base point
$p_k$.
This determines
the fundamental groupoid
$\widetilde \pi= \Pi(\Sigma;p_0,p_1,\dots,p_n)$.
To obtain a presentation of it we decompose
$\Sigma$ into cells as follows
where we do not distinguish in notation between the chosen
edge paths and their homotopy classes relative
to their end points:
Let $x_1,y_1,\dots,x_{\ell}, y_{\ell}$
be closed paths 
which (i) do not meet the boundary,
(ii) have $p_0$ as starting point, and
(iii) yield the generators
respectively
$x_1,y_1,  \dots,x_\ell,y_\ell$
of the fundamental group
$\pi=\pi_1(\Sigma,p_0)$;
for $j = 1,\dots, n$,
let $a_k$ 
be
the boundary path of the $k$'th
boundary circle,
having $p_k$ as starting point,
and let
$\gamma_k$ 
be a path
from $p_0$ to $p_k$.
When we cut $\Sigma$ along these
1-cells we obtain a disk whose boundary
yields the defining relation of 
$\widetilde \pi= \Pi(\Sigma;p_0,p_1,\dots,p_n)$.
The resulting
presentation
of 
$\widetilde \pi$
looks like
$$
\widetilde {\Cal P} = \langle x_1,y_1,  \dots,x_\ell,y_\ell,a_1,\dots,a_n,
\gamma_1,\dots,\gamma_n; \widetilde r\rangle,
\quad \widetilde r 
= \Pi_{j=1}^{\ell} [x_j,y_j] \Pi_{k=1}^n \gamma_k a_k \gamma_k^{-1}.
\tag2.1
$$
Let
$\widetilde F$ denote the  groupoid
which is free
on the generators of (2.1).
To have a neutral notation,
whenever necessary,
we shall write
$\widetilde \Pi$
for either
$\widetilde F$ or
$\widetilde \pi$;
accordingly we write
$\Pi$
for either
$F$ or
$\pi$.
This is consistent with the presentation (1.1) of
the fundamental group if we identify,
for $1 \leq k \leq n$, the generator
$z_k$ with $\gamma_k a_k \gamma_k^{-1}$.
More precisely,
the assignments
$$
\aligned
i(e) &= p_0,
\quad
i(x_j) = x_j,
\quad
i(y_j) = y_j,
\quad
i(z_k) = \gamma_k a_k \gamma_k^{-1},
\\
\rho(p_j) &= e,
\quad
\rho(x_j) = x_j,
\quad
\rho(y_j) = y_j,
\quad
\rho(a_k) = z_k,
\quad
\rho(\gamma_k) =\roman{Id},
\endaligned
\tag2.2
$$
where $1 \leq j\leq \ell$ and $1 \leq k \leq n$,
yield 
morphisms of presentations
$i \colon \Cal P \to \widetilde {\Cal P}$
and
$\rho \colon \widetilde {\Cal P} \to \Cal P$
and hence
functors 
$i \colon \Pi \to \widetilde \Pi$
and
$\rho \colon \widetilde \Pi \to \Pi$
(where the notation $i$ and $\rho$ is abused)
inducing a deformation retraction
of $\widetilde \Pi$ onto $\Pi$;
cf. e.~g. \cite{\brownboo\ (6.5.13)}
for the latter notion.
\smallskip
As usual, view the group
$G$ as a groupoid with a single object $e$, 
identified with the neutral element of $G$.
Denote by
$\widetilde \Pi_0$ the set
$\{p_0,\dots,p_n\}$
of objects
of $\widetilde \Pi$, and
write $\roman{Hom}(\widetilde \Pi,G)$
for the space of groupoid homomorphisms from
$\widetilde \Pi$ to $G$.
The obvious action of $G$
on $\roman{Hom}(\Pi,G)$
by conjugation
extends
to an action of
the group
$G^{\widetilde \Pi_0} \cong G\times \dots \times G$ 
($n+1$ copies of $G$)
on $\roman{Hom}(\widetilde \Pi,G)$
in the following way:
We denote by $s$ and $t$ the source and target mappings
from $\widetilde \Pi$ to 
the object set
$\widetilde \Pi_0$.
Given a homomorphism
$\alpha$ from
$\widetilde \Pi$ to $G$ and 
$\vartheta \in G^{\widetilde \Pi_0}$,
the homomorphism
$\vartheta \alpha$ is defined by
$$
\vartheta \alpha(w) = \vartheta(t(w))\alpha(w)(\vartheta(s(w)))^{-1}.
$$
The orbit space for the 
$G^{\widetilde \Pi_0}$-action on
$\roman{Hom}(\widetilde \Pi,G)$
will be denoted by 
$\roman{Rep}(\widetilde \Pi,G)$.
The structure of 
the space $\roman{Rep} (\Pi,G)$ 
of representations of the group $\Pi$
can now be studied
by looking at
$\roman{Rep}(\widetilde \Pi,G)$
instead. More precisely:
The  functors $i$ and $\rho$ induce 
maps
$$
i^* \colon \roman{Hom}(\widetilde \Pi,G)
@>>>
\roman{Hom} (\Pi,G),
\quad
\rho^* \colon \roman{Hom}( \Pi,G)
@>>>
\roman{Hom} (\widetilde \Pi,G)
$$
which, for $\Pi = F$ and
$\widetilde \Pi = \widetilde F$, are manifestly smooth.
We shall occasionally refer to
$i^*$ and
$\rho^*$ as
{\it restriction\/}
and
{\it corestriction\/},
respectively.
The following is obvious.

\proclaim{Proposition 2.3}
The restriction 
mapping
induces a homeomorphism
$$
i^* \colon \roman{Rep}(\widetilde \Pi,G)
@>>>
\roman{Rep} (\Pi,G). \qed
$$
\endproclaim

\smallskip
We apply a variant of 
the construction in \cite\modus\ to the presentation
$\widetilde {\Cal P}$:
Let $O_0, O_1,\dots, O_n$ be open 
$G$-invariant subsets of the Lie algebra
$\fra g$ of $G$ where the exponential map is regular;
at this stage,
we just take
for
each of 
$O_0, O_1,\dots, O_n$ the {\it same\/} maximal open connected
$G$-invariant neighborhood of the origin
of $\fra g$
where the exponential map is regular;
the subscripts then constitute merely a notational device.
In (2.18) below we shall take
$O_1 =\dots = O_n = O$,
where $O$ refers to the chosen
invariant connected neighborhood of zero
where the exponential map
is a diffeomorphism.
With the present
choice
of $O_0, O_1,\dots, O_n$,
define
the space
$\Cal H(\widetilde {\Cal P},G)$
by means of the pull back square
$$
\CD
\Cal H(\widetilde {\Cal P},G)
@>{(\widehat {\widetilde r},\widehat a_1,\dots,\widehat a_n)}>>
O_0 \times O_1\times\dots\times O_n
\\
@V{{\widetilde \eta}}VV
@VV{\roman{exp} \times \dots\times \roman{exp}}V
\\
\roman{Hom}(\widetilde F,G)
@>>{(\widetilde r,a_1,\dots,a_n)}> G_0\times G_1 \times \dots \times G_n,
\endCD
\tag2.4
$$
where
$\widehat {\widetilde r}$ and $\widehat a_1,\dots,\widehat a_n$
denote the smooth maps
induced by, respectively,
$\widetilde r$ and $a_1,\dots,a_n$,
and where we have written
$G_0=G, \, G_1=G, \dots, G_n = G$;
here $\widetilde \eta$ is just a name for the corresponding map
which results from the pull back construction.
Since $\roman{Hom}(\widetilde F,G)$
is a smooth manifold, and since
the  exponential map, restricted to any of the $O_j$'s, is regular,
$1 \leq j \leq n$,
the space $\Cal H(\widetilde {\Cal P},G)$
is a smooth manifold, too.
\smallskip
Keeping the notation used in \cite\modus,
for a group or groupoid $\Pi$, we denote by $(C_*(\Pi,R),\partial)$
the chain complex of its inhomogeneous reduced normalized
bar resolution over $R$.
Let $c$ be an absolute 2-chain of $F$ 
which represents a 2-cycle for the group system
$(\pi;\pi_1,\dots,\pi_n)$. 
Its image in the 2-chains
of the fundamental group
$\widehat \pi$
of the {\it closed\/} (!)
surface $\widehat \Sigma$ is then closed.
Write $\widehat\kappa \in \roman H_2(\widehat \pi)$ for its class.
When the genus $\ell$ is different from zero,
$\widehat \pi$ is non-trivial and
the canonical map
from
the second homology group
$\roman H_2(\pi,\{\pi_k\})$
of the group system (cf. \cite\guhujewe)
to $\roman H_2(\widehat \pi)$
is an isomorphism
identifying the fundamental classes.
Furthermore, when the relators 
$z_1,\dots,z_n$ are added to
(1.1) we obtain the presentation
$$
\widehat {\Cal P} = \langle x_1,y_1,  \dots,x_\ell,y_\ell,z_1,\dots,z_n; r,
z_1,\dots,z_n
\rangle
\tag2.5
$$
of the fundamental group
$\widehat \pi =\pi_1(\widehat \Sigma)$
of the {\it closed\/}
surface $\widehat \Sigma$ resulting from capping of the $n$ boundaries.
\smallskip
We now apply a variant of the construction
in Theorem 1 
of \cite \modus:
Write
$\widetilde \pi_{\partial}$
for the free subgroupoid
of $\widetilde F$
having $p_1,\dots,p_n$ as objects and
$a_1,\dots,a_n$ as morphisms;
this groupoid may also be viewed as a subgroupoid of
the fundamental groupoid
$\widetilde \pi$ of $\Sigma$, 
and we do no distinguish in notation between the two 
subgroupoids.
Abstractly,
$\widetilde \pi_{\partial}$
amounts of course to a disjoint union
of the $n$ free cyclic groups
$\pi_1,\dots, \pi_n$.
We pick
$c$ in such a way that
$$
\partial c = [r] - [z_1] - \dots - [z_n]
\tag2.6
$$
in the reduced normalized
inhomogeneous
bar resolution of $F$.
This can always be done,
cf. what is said 
in Section 2 of \cite\guhujewe.
View $c$ as a 2-chain of
$\widetilde F$
by the embedding $\iota$ of
$F$ into $\widetilde F$
(cf. (2.2))
and let
$$
\widetilde c =
c + \sum_j\left([\gamma_j^{-1}|\gamma_j a_j] - 
[\gamma_ja_j|\gamma_j^{-1}]
\right),
\tag2.7
$$
cf. \cite{\guhujewe\  (8.7)}.
Then
$$
\partial \widetilde c =
[\widetilde r] - [a_1] - \dots - [a_n]
\tag2.8
$$
whence, in particular, $\widetilde c$ is 
manifestly
a relative cycle
for $(\widetilde \pi, \widetilde \pi_{\partial})$,
cf. step 3 in the proof of the key lemma
(8.4)
in \cite\guhujewe.
Notice that
$c$ itself is {\it not\/} 
a relative cycle
for $(\widetilde \pi, \widetilde \pi_{\partial})$.
\smallskip
Let $\cdot$ be an invariant symmetric bilinear form on $\fra g$,
at this stage not necessarily nondegenerate,
and let $\Omega = \omega_1 \cdot \overline{\omega}_2\in \Omega^2 (G^2)$
be the indicated 2-form 
arising from the Maurer-Cartan form $\omega$ on $G$ with values
in $\fra g$,
cf. \cite{\guhujewe,\ \modus,\ \weinsthi}.
Write 
$$
E\colon \widetilde F^2 \times \roman{Hom}(\widetilde F,G) @>>> G^2
$$
for the evaluation map, and let
$$
\omega_{\widetilde c} = \langle \widetilde c, E^*\Omega\rangle,
\tag2.9
$$
the result of pairing $\widetilde c$ with the induced form,
cf. 
\cite{\guhujewe\  (5.3)}
and
\cite {\modus\  (13)}.
This is a $G$-invariant  2-form on
$\roman{Hom}(\widetilde F,G)$.
In view of 
\cite {\modus\  (15)}
we have
$$
d\omega_{\widetilde c} = \langle \partial \widetilde c, E^*\lambda\rangle.
\tag2.10
$$
Here $\lambda$ is the fundamental 3-form on $G$,
cf. \cite {\guhujewe,\ \modus,\ \weinsthi}.
Let 
$h$ be the standard homotopy operator
on forms on $\fra g$ arising from integration along straight
line segments,
let $\beta = h(\roman{exp}^* \lambda)$, and 
define
the 2-form $\omega_{\widetilde c,\widetilde {\Cal P}}$ on
$\Cal H(\widetilde {\Cal P},G)$ by
$$
\omega_{\widetilde c,\widetilde {\Cal P}}=
{\widetilde \eta}^* \omega_{\widetilde c}
-
\widehat {\widetilde r}^* \beta
+
\widehat a_1^* \beta
+
\dots
+
\widehat a_n^* \beta.
\tag2.11.1
$$
In view of (2.8) and (2.10),
this 2-form is closed;
cf. the reasoning in Section 7 of \cite\guhujewe.
\smallskip
Let $\psi \colon \fra g \to \fra g^*$
be the adjoint of
the given
invariant symmetric bilinear form $\cdot$ on $\fra g$,
and write 
$$
\mu \colon
\Cal H(\widetilde {\Cal P},G)
@>>>
\fra g_0^* \times \fra g^*_1\times\dots\times \fra g^*_n
\tag2.11.2
$$
for the composite
$$
\Cal H(\widetilde {\Cal P},G)
@>{(\widehat {\widetilde r},\widehat a_1,\dots,\widehat a_n)}>>
O_0 \times \dots\times O_n
@>>>
\fra g_0 \times\dots\times \fra g_n
@>{\psi^{\times (n+1)}}>>
\fra g_0^* \times\dots\times \fra g^*_n.
$$
Further,
for $0 \leq j \leq n$, we write
$$
\mu_j \colon
\Cal H(\widetilde {\Cal P},G)
@>>>
\fra g_j^*
\tag 2.11.3
$$
for the composite
of $\mu$ with the projection onto the $j$'th copy 
$\fra g_j^*$
of
$\fra g^*$
in
$\fra g_0^* \times\dots\times \fra g^*_n$.
The same reasoning as in Section 7 of \cite\guhujewe\ 
shows that
$$
-\omega_{\widetilde c,\widetilde {\Cal P}}
(X_{\Cal H},\cdot)
= d (X \circ \mu),
\tag2.12.1
$$
that is to say,
the formal momentum mapping property is satisfied
where we have written $\Cal H =\Cal H(\widetilde {\Cal P},G)$
for short.
Here
$X_{\Cal H}$
denotes the vector field
on $\Cal H$ coming from
$X \in \fra g_0 \times \dots \times \fra g_n$ 
via the 
action
of $G_0 \times \dots \times G_n$
on $\Cal H$.
Likewise,
for $0 \leq j \leq n$,
$\mu_j$ is formally a momentum mapping for the corresponding
action of 
the $j$'th copy 
$G_j$ of $G$ on
$\Cal H(\widetilde {\Cal P},G)$,
that is,
$$
-\omega_{\widetilde c,\widetilde {\Cal P}}
(X_{\Cal H},\cdot)
= d (X \circ \mu_j),
\tag2.12.2
$$
where now
$X_{\Cal H}$
refers to the vector field
on $\Cal H$ coming from
$X \in \fra g_j$
via the 
$G_j$-action
on $\Cal H$.
\smallskip
We now suppose that the 
given
invariant symmetric bilinear form $\cdot$ on $\fra g$
is nondegenerate;
this allows us to identify $\fra g$ tacitly with its dual $\fra g^*$.
In particular, 
this will always be meant 
below
when we refer to the \lq\lq preimage of the center 
with respect to the momentum mapping\rq\rq\ 
and when we identify adjoint orbits with coadjoint ones.

\proclaim{Theorem 2.13}
Near the zero locus of the map
$\widehat{\widetilde r}\colon \Cal H(\widetilde {\Cal P},G) \to O_0$
brought into play in {\rm (2.4)}
or, more generally,
near the preimage 
with respect to $\widehat{\widetilde r}$
of the center of $\fra g$,
the 2-form
$\omega_{\widetilde c,\widetilde {\Cal P}}$
on
$\Cal H(\widetilde {\Cal P},G)$,
cf. {\rm (2.11.1)},
is nondegenerate, that is, 
a symplectic structure.
\endproclaim

We note that, 
when $O_1 =\dots = O_n = O$,
the chosen
invariant connected neighborhood of zero, $O$,
where the exponential map is a diffeomorphism,
the smooth map $\widetilde \eta$ (cf. (2.4))
is injective and
the zero locus of $\widehat{\widetilde r}$
may then be identified
with the space $\roman{Hom}(\widetilde \pi, G)$.
\smallskip
To prepare for the proof, which will be given in (2.17) below,
and to gain
additional insight,
we proceed as follows:
Let $\Cal O= (\Cal O_1,\dots, \Cal O_n)$
be an $n$-tuple of (co)adjoint orbits, each
$\Cal O_k$ being in $O_k$,
and let
$\bold C= (C_1,\dots, C_n)$
be the corresponding $n$-tuple of conjugacy classes,
where $C_k = \roman{exp}(\Cal O_k)$,
for $ 1 \leq k \leq n$.
We may in fact identify
$\Cal O$ with 
the (co)adjoint orbit
$\Cal O_1 \times \dots \times \Cal O_n$
for
$\Gamma = G_1 \times \dots \times G_n$
and
$\bold C$
with the 
conjugacy class
$C_1 \times \dots \times C_n$
of
$\Gamma$.
Since 
each $\Cal O_k$ is in $O_k$,
each projection map
$\Cal O_k \to C_k$ is a covering projection.
Consider the mapping
$$
\mu_1 \times \mu_2\times \dots \times \mu_n 
\colon
\Cal H(\widetilde {\Cal P},G)
@>>>
\fra g^*_1\times\dots\times \fra g^*_n;
$$
in view of (2.12.1) and (2.12.2),
this is a momentum mapping
for the action of
$\Gamma$
on $\Cal H(\widetilde {\Cal P},G)$,
with reference to the closed 2-form
$\omega_{\widetilde c,\widetilde {\Cal P}}$.
Whether or not the latter is nondegenerate, we can form the reduced space
$$
\Cal H(\widetilde {\Cal P},G)_{\overline{\Cal O}}
=
\Cal H(\widetilde {\Cal P},G)_
{\overline{\Cal O}_1 \times \dots \times\overline{\Cal O}_n};
$$
with reference to the $n$-tuple
$$
\omega = (\omega_1,\dots, \omega_n)
$$
of the Kirillov forms
$\omega_k$ on the $\Cal O_k$'s,
the reduced space $\Cal H(\widetilde {\Cal P},G)_{\overline{\Cal O}}$ 
inherits a closed 2-form 
$$
\omega_{\widetilde c,\widetilde {\Cal P},\omega},
\tag2.14.1
$$
the corresponding {\it reduced\/} form, in the standard way.
As usual,  the notation
$\overline{\Cal O}_1$,
$\overline{\Cal O}_2$ etc. here indicates
the (co)adjoint orbit
$\Cal O_1$,
$\Cal O_2$ etc., 
endowed with the {\it negative\/}
of the symplectic (i.~e. Kirillov)
structure of
$\Cal O_1$, $\Cal O_2$ etc.
We may in fact view
$\omega$
as the Kirillov form of the (co)adjoint orbit
$\Cal O_1 \times \dots \times \Cal O_n$
for
$\Gamma = G_1 \times \dots \times G_n$.
The hamiltonian action of
$G_0 \times \dots\times G_n$ 
on $\Cal H(\widetilde {\Cal P},G)$
induces a hamiltonian
action
of $G_0$ on
$\Cal H(\widetilde {\Cal P},G)_{\overline{\Cal O}}$,
preserving
$\omega_{\widetilde c,\widetilde {\Cal P},\omega}$
and
having momentum mapping
$$
\mu_{\Cal O}
\colon
\Cal H(\widetilde {\Cal P},G)_{\overline{\Cal O}}
@>>>
\fra g_0^*
\tag2.14.2
$$
induced by $\mu_0$ or, what amounts to the same,
induced by the map
$\widehat{\widetilde r} \colon
\Cal H(\widetilde {\Cal P},G)
\to O_0$ given in (2.4).
As a space,
$$
\Cal H(\widetilde {\Cal P},G)_{\overline{\Cal O}}
=
\left(\mu_1 \times  \dots \times \mu_n\right)^{-1}
\left(\Cal O_1 \times  \dots \times \Cal O_n\right)\big /
(G_1\times \dots \times G_n).
$$
On the other hand, 
define the smooth manifold
$\Cal H({\Cal P},G)_{\overline{\Cal O}}$
by means of the pull back square
$$
\CD
\Cal H({\Cal P},G)_{\overline{\Cal O}}
@>{(\widehat r,\widehat z_1,\dots,\widehat z_n)}>>
O_0 \times \Cal O_1\times\dots\times \Cal O_n
\\
@V{\eta^{\flat}}VV
@VV{\roman{exp} \times \roman{exp} \times \dots\times \roman{exp}}V
\\
\roman{Hom}(F,G)_{\bold C}
@>>{(r,z_1,\dots,z_n)}> G_0\times C_1 \times \dots \times C_n,
\endCD
\tag2.14.3
$$
where the notation
$\widehat r,\widehat z_1,\dots,\widehat z_n$ etc. is the obvious one.
When the action of
$\Gamma$ is divided out,
the groupoid generators
$\gamma_k$ are no longer relevant
in the sense that the functor $i$
from $F$ to $\widetilde F$
spelled out in  (2.2)
induces a diffeomorphism
from
$\Cal H(\widetilde {\Cal P},G)_{\overline{\Cal O}}$
onto
$\Cal H({\Cal P},G)_{\overline{\Cal O}}$,
cf. also (2.3),
and we identify henceforth the two spaces;
notice in particular that at this stage
the distinction 
in notation between the $a_k$'s and the $z_k$'s
is no longer necessary.
\smallskip
In \cite{\guhujewe\ (5.2)},
for arbitrary conjugacy classes
$C_1,\dots,C_n$,
(i.e. not necessarily consisting of regular points
for the exponential map),
the space
$\Cal H(\wide {\Cal P},G)_{\bold C}$,
from which the extended moduli in that paper 
is obtained,
is defined by means of the pull back square
$$
\CD
\Cal H(\wide {\Cal P},G)_{\bold C}
@>{(\widehat r,\overline z_1,\dots,\overline z_n)}>>
O_0 \times C_1\times\dots\times C_n
\\
@V{\eta}VV
@VV{\roman{exp} \times \roman{Id} \times \dots\times \roman{Id}}V
\\
\roman{Hom}(F,G)_{\bold C}
@>>{(r,z_1,\dots,z_n)}> G_0\times C_1 \times \dots \times C_n,
\endCD
\tag2.14.4
$$
where
$\widehat r$ and $\overline z_1,\dots,\overline z_n$
denote the smooth maps
induced by,
respectively,
$r$ and $z_1,\dots,z_n$;
again  $\eta$ is just a name for the corresponding map
which results from the pull back construction.
Hence
the obvious morphism of squares from
(2.14.3) to (2.14.4)
and 
the functor $i$
from $F$ to $\widetilde F$
plainly induce
a smooth mapping 
$$
i^{\sharp}\colon
\Cal H(\widetilde {\Cal P},G)_{\overline{\Cal O}}
@>>>
\Cal H(\wide {\Cal P},G)_{\bold C}
\tag2.14.5
$$
of $G$-spaces, $G$ being, roughly speaking, the residual copy $G_0$ of $G$
in $G_0\times G_1 \times \dots \times G_n$.
Recall that the corresponding momentum mapping
$$
\mu_{\bold C}
\colon
\Cal H({\Cal P},G)_{\bold C}
@>>>
\fra g_0^*
\tag2.14.6
$$
is induced by 
$\widehat r \colon \Cal H({\Cal P},G)_{\bold C} \to O_0$;
in \cite{\guhujewe\  (Section 7)}
this momentum mapping is written
$\mu\colon\Cal H({\Cal P},G)_{\bold C}
\to \fra g^*$.
It is obvious that 
$i^{\sharp}$ is compatible with
the momentum mappings
(2.14.2) and (2.14.6)
since both come from the relator $r$.
We remind the reader that, at this stage,
for each $j$, the exponential map, restricted
to $\Cal O_j$, is a covering projection onto $C_j$,
for $0 \leq j \leq n$.

\proclaim{Theorem 2.14}
Near the zero locus of 
the momentum  mapping
$\mu_{\Cal O}$ brought into play in
{\rm (2.14.2)}
or, more generally,
near the preimage 
with respect to $\mu_{\Cal O}$
of the center of $\fra g$,
the reduced 2-form 
$\omega_{\widetilde c,\widetilde {\Cal P},\omega}$
on
$\Cal H(\widetilde {\Cal P},G)_{\overline{\Cal O}}$,
cf. {\rm (2.14.1)},
is symplectic, and
the map
$i^{\sharp}$
from 
$\Cal H(\widetilde {\Cal P},G)_{\overline{\Cal O}}$
to
$\Cal H({\Cal P},G)_{\bold C}$
is a covering projection of symplectic manifolds
with hamiltonian $G$-action.
\endproclaim

The proof will be postponed until (2.17) below;
some prerequisites for it will be given in (2.15) and (2.16).

\smallskip\noindent
(2.15) 
{\smc Presymplectic reduction\/}.
Let $M$ be a smooth manifold endowed with a 
presymplectic
(i.~e. closed)  2-form $\omega$
and an action  of a Lie group $K$ which is hamiltonian in the sense
that there is a smooth $K$-equivariant map
$\mu$ from $M$ to the dual  $\fra k^*$
of the Lie algebra $\fra k$
which satisfies
the formal momentum property
$$
-\omega(X_M, \cdot)
= d (X \circ \mu) (= X \circ d \mu),
$$
where $X_M$ refers to the vector field
on $M$ coming from
$X \in \fra k$
via the $K$-action.
Let $p$ be a point of $M$, 
and let
$$
Z_p = \roman{ker}(d\mu_p) \subseteq \roman T_pM,
\quad
B_p = \roman T_p(Gp) \subseteq \roman T_pM.
$$
As usual, we write $B_p^\omega$ etc. for the 
{\it annihilator\/}
of $B_p$ in $\roman T_p M$, with reference to $\omega$.

\proclaim{2.15.1}
The annihilator of
$B_p$ coincides with $Z_p$, that is,
$B_p^\omega = Z_p$.
\endproclaim

\noindent
Next we suppose that $\mu(p) = 0$, that is to say,
that $p$ lies in the zero locus of the momentum mapping $\mu$.

\proclaim{2.15.2}
The vector space $B_p$ is isotropic,
that is, 
$B_p \subseteq  B_p^\omega = Z_p$.
\endproclaim

\noindent
In fact, the restriction of $\mu$ to the orbit $Gp$
is constant, having constant value zero. Hence
the differential $d\mu_p$,
restricted to $B_p$ is zero whence,
for $X,Y \in \fra k$,
$$
-\omega(X_M,Y_M)_p = X d\mu_p(Y_M) = 0.
$$

\proclaim{2.15.3}
The vector space $Z_p$ is coisotropic,
that is, 
$Z_p^\omega \subseteq Z_p$.
\endproclaim

\noindent
In fact, since 
$B_p \subseteq  Z_p$, we have
$Z_p^\omega \subseteq  B_p^\omega = Z_p$.---Let $H_p = Z_p /B_p$.

\proclaim{Lemma 2.15.4}
The 2-form  $\omega$ is nondegenerate near the point $p$ 
in the zero locus of $\mu$
if and only
if the 
the following holds:
The alternating bilinear form
$\omega_p$
on $H_p$ induced by $\omega$
is nondegenerate and
the image $d\mu_p(\roman T_pM)$ in $\fra k^*$
consists precisely of the linear forms on $\fra k$
which annihilate the Lie algebra $\fra k_p$ of the stabilizer
$K_p$ of the point $p$ of $M$.
\endproclaim

\demo{Proof} 
It is a standard fact that the condition is necessary.
Sufficiency is readily verified using a variant of the argument
given in Section 5 of \cite\modus:
When $\omega_p$ is symplectic,
the alternating
bilinear form $\omega$ on $\roman T_pM$, restricted to $Z_p$, has degeneracy
space equal to the subspace $B_p$ of $Z_p$.
The momentum mapping property,
combined with the condition involving the stabilizer
$\fra k_p$ at $p$,
then implies that
on
$\roman T_pM$,
the 
form $\omega$ 
is nondegenerate:
Let $X \in \roman T_pM$, and suppose that
$$
\omega(X,Y) = 0,
\quad
\text{for every}
\quad
Y \in \roman T_pM.
$$
Taking $Y=U_M$ for $U \in \fra k$ and
using the formal momentum mapping property,
we see that
$U d\mu_p (X) = 0$
for every $U \in \fra k$
whence
$d\mu_p (X) = 0$
and thence
$X \in Z_p$.
The nondegeneracy of $\omega_p$ on $H_p$
then implies that $X \in B_p$,
that is, $X = V_M$ for some $V \in \fra k$.
The formal momentum mapping property
then yields
$V d\mu_p (Y) = 0$
for every
$Y \in \roman T_pM$,
that is, $V$ annihilates the image
$d\mu_p(\roman T_pM)$ in $\fra k^*$.
Hence
$V$ lies in $\fra k_p$
whence
$X=V_M$ vanishes at the point $p$. \qed
\enddemo

\smallskip\noindent
(2.16) {\smc Parabolic cohomology}. 
Let $\phi \in \roman{Hom}(F,G)$,
and suppose that $\phi(r)$ lies in the center of $G$.
Then the composite
of $\phi$
with the adjoint action
of $G$ induces
the structure
of a left $\pi$-module 
on the Lie algebra $\fra g$, and we write
$\fra g_{\phi}$
for $\fra g$, viewed as a $\pi$-module in this way.
The Reidemeister-Fox
calculus, applied to
the presentation (1.1), 
yields  free resolution
$\bold R(\Cal P)$,
cf. \cite{\guhujewe\ (2.2)},
and application of the functor
$\roman{Hom}_{\Bobb R\pi}(\cdot, \fra g_{\phi})$
to this free resolution
yields the chain complex
$$
\bold C(\Cal P, \fra g_{\phi})\colon
\roman C^0(\Cal P, \fra g_{\phi})
@>{\delta_{\phi}^0}>>
\roman C^1(\Cal P, \fra g_{\phi})
@>{\delta_{\phi}^1}>>
\roman C^2(\Cal P, \fra g_{\phi}),
\tag2.16.1
$$
cf. \cite{\modus\ (4.1)}
and \cite{\guhujewe\  (4.1)};
this chain complex 
computes the group cohomology $\roman H^*(\pi,\fra g_{\phi})$. 
We note that
there are canonical isomorphisms
$$
\roman C^0(\Cal P, \fra g_{\phi}) \cong \fra g,
\quad
\roman C^1(\Cal P, \fra g_{\phi}) \cong \fra g^{2\ell+ n},
\quad
\roman C^2(\Cal P, \fra g_{\phi}) \cong \fra g .
$$
To recall the geometric significance
of this chain complex,
denote by $\alpha_\phi$
the smooth map 
from $G$  to $\roman{Hom}(F,G)$
which assigns $x \phi x^{-1}$ to $x \in G$,
maintain the notation $r$
for the smooth map
from
$\roman{Hom}(F,G)$
to
$G$
induced by the relator $r$
so that the pre-image of the neutral element
$e$ of $G$ equals
the space
$\roman{Hom}(\pi,G)$,
and write
$R_\phi\colon \fra g^{2\ell+n} \to  \roman T_{\phi} \roman{Hom}(F,G)$
and
$R_{r\phi}\colon \fra g \to  \roman T_{r (\phi)}G$
for the corresponding operations of right translation.
The tangent maps
$\roman T_e\alpha_{\phi}$
and $\roman T_{\phi}r$
make commutative the diagram
$$
\CD
\roman T_eG
@>\roman T_e\alpha_{\phi}>>
\roman T_{\phi} \roman{Hom}(F,G)
@>{\roman T_{\phi} r}>>
\roman T_{r (\phi)}G
\\
@A{\roman{Id}}AA
@A{\roman R_{\phi}}AA
@A{\roman R_{r(\phi)}}AA
\\
\fra g
@>>{\delta^0_{\phi}}>
\fra g^{2\ell +n}
@>>{\delta^1_{\phi}}>
\fra g,
\endCD
\tag2.16.2
$$
cf. \cite{\guhujewe\ (4.2)} and \cite{\modus\  (4.2)}.
We now suppose that   our chosen
$\phi \in \roman{Hom}(F,G)$ lies in 
$\roman{Hom}(F,G)_{\bold C}$,
viewed as a subspace
of
$\roman{Hom}(F,G)$.
For $k = 1, \dots, n$,
denote by $\fra h_k$ the image in $\fra g$ of the linear endomorphism
given by
$\roman{Ad}(\phi(z_k)) -\roman{Id}$, so that there results the exact sequence
$$
0
@>>>
\fra s_k
@>>>
\fra g
@>>>
\fra h_k
@>>>
0
\tag2.16.3
$$
of vector spaces,
where
$\fra s_k$ denotes the 
Lie algebra of the
stabilizer
of $\phi(z_k)$;
notice that $\fra h_k$ amounts to the tangent space
of the conjugacy class $C_k$.
By
Proposition 4.4
of \cite\guhujewe,
the values of the operator
$\delta^0_{\phi}$
in (2.16.2)
lie in
$\fra g^{2\ell} \times \fra h_1 \times \dots\times \fra h_n$, 
viewed as a subspace
of
$\roman C^1(\Cal P, \fra g_{\phi}) \cong \fra g^{2\ell} \times \fra g^n$,
and
the first cohomology group of the resulting chain complex
$$
\bold C_{\roman{par}}(\Cal P, \fra g_{\phi})\colon
\fra g
@>{\delta^0_{\phi}}>>
\fra g^{2\ell} \times \fra h_1 \times \dots\times \fra h_n
@>{\delta^1_{\phi}}>>
\fra g
\tag2.16.4
$$
equals 
the first parabolic cohomology group
$\roman H_{\roman{par}}^1(\pi,\{\pi_j\};\fra g_{\phi})$
in the sense of \cite\weiltwo;
the latter may in fact be defined to be
the image of the relative cohomology group
$\roman H^1(\pi,\{\pi_j\};\fra g_{\phi})$
in the absolute cohomology group
$\roman H^1(\pi,\fra g_{\phi})$ under the canonical map.
We now extend this notion
by defining the {\it parabolic cohomology\/}
$\roman H_{\roman{par}}^*(\pi,\{\pi_j\};\fra g_{\phi})$
of
$(\pi,\{\pi_j\})$ {\it with values in\/} $\fra g_{\phi}$
to be the cohomology of (2.16.4);
in view of the cited Proposition 4.4 in \cite\guhujewe,
this yields the correct notion in degree 1.
Further, it is straightforward to generalize this definition
of parabolic cohomology
to an arbitrary $\pi$-module $V$ but we shall not need this
here.
\smallskip
We now determine 
$\roman H_{\roman{par}}^0(\pi,\{\pi_j\};\fra g_{\phi})$
and 
$\roman H_{\roman{par}}^2(\pi,\{\pi_j\};\fra g_{\phi})$.
Let
$$
\Phi=
(\roman{Ad}(z_1)-\roman{Id},\dots,\roman{Ad}(z_n)-\roman{Id})
\colon
\fra g^n
@>>>
\fra g^n.
$$
Application of the functor $\roman {Hom}_{\bold R \pi}(\cdot,\fra g_{\phi})$
to the 
quotient complex
$\bold R(\widetilde{\Cal P}, \{\pi_j\})$
constructed in \cite{\guhujewe\ (2.8)}
yields the cochain complex
$$
\bold C(\widetilde {\Cal P},\{\pi_j\}; \fra g_{\phi}) \colon
\fra g
@>{\delta_{\phi}^0}>>
\fra g^{2\ell}  \times \fra g^n
@>{\delta_{\phi}^1}>>
\fra g,
$$
and application of this functor to
the comparison map
\cite{\guhujewe\ (2.12)} yields the cochain map
$$
(\roman{Id}, (\roman {Id}, \Phi),\roman{Id})
\colon
\bold C(\widetilde {\Cal P},\{\pi_j\}; \fra g_{\phi})
@>>>
\bold C(\Cal P, \fra g_{\phi})
\tag2.16.5
$$
which induces the canonical map
$$
\roman H^*(\pi,\{\pi_j\}; \fra g_{\phi})
@>>>
\roman H^*(\pi, \fra g_{\phi}).
$$
By construction,
(2.16.5) may be factored through
(2.16.4); when we display this, we obtain the commutative diagram
$$
\CD
\bold C(\Cal P, \fra g_{\phi}) \colon
@.
\fra g
@>{\delta_{\phi}^0}>>
\fra g^{2\ell}  \times \fra g^n
@>{\delta_{\phi}^1}>>
\fra g
\\
@.
@A{\roman{Id}}AA
@A{(\roman{Id},\roman{incl})}AA
@A{\roman{Id}}AA
\\
\bold C_{\roman{par}}(\Cal P, \fra g_{\phi})\colon
@.
\fra g
@>{\delta^0_{\phi}}>>
\fra g^{2\ell} \times \fra h_1 \times \dots\times \fra h_n
@>{\delta^1_{\phi}}>>
\fra g
\\
@.
@A{\roman{Id}}AA
@A{(\roman{Id},\roman{proj})}AA
@A{\roman{Id}}AA
\\
\bold C(\widetilde {\Cal P},\{\pi_j\}; \fra g_{\phi}) \colon
\quad
@.
\fra g
@>>{\delta_{\phi}^0}>
\fra g^{2\ell}  \times \fra g^n
@>>{\delta_{\phi}^1}>
\fra g
\endCD
\tag2.16.6
$$
where the notation 
$\delta_{\phi}^0$
and
$\delta_{\phi}^1$
is slightly abused;
here \lq\lq proj\rq\rq\ and
\lq\lq incl\rq\rq\ refer to the
obvious projection and inclusion mappings,
and the composite
$\roman{incl} \circ \roman{proj}$
equals $\Phi$.
From the commutativity of
(2.16.6)
we deduce at once the following canonical isomorphisms
$$
\roman H_{\roman{par}}^0(\pi,\{\pi_j\};\fra g_{\phi})
\cong
\roman H^0(\pi,\fra g_{\phi});
\quad
\roman H_{\roman{par}}^2(\pi,\{\pi_j\};\fra g_{\phi})
\cong
\roman H^2(\pi,\{\pi_j\};\fra g_{\phi}).
\tag2.16.7
$$
Poincar\'e duality for the surface group system
$(\pi; \pi_1,\dots, \pi_n)$
yields the duality isomorphism
$$
\cap \kappa
\colon
\roman H^2(\pi,\{\pi_j\};\fra g_{\phi})
@>>>
\roman H_0(\pi,\fra g_{\phi})
$$
where $\kappa$ refers to the fundamental class in
$\roman H_2(\pi,\{\pi_j\};\Bobb R)$,
cf. Section 3 of \cite\guhujewe.
Consequently Poincar\'e duality
yields a canonical duality isomorphism
$$
\cap \kappa
\colon
\roman H_{\roman{par}}^2(\pi,\{\pi_j\};\fra g_{\phi})
@>>>
\roman H_0(\pi,\fra g_{\phi}).
\tag2.16.8
$$
Since the chosen invariant symmetric nondegenerate bilinear from
$\cdot$ on $\fra g$ induces 
a nondegenerate
bilinear
pairing
$$
\roman H^0(\pi,\fra g_{\phi})
\otimes
\roman H_0(\pi,\fra g_{\phi})
@>>>
\Bobb R
$$
as usual, from (2.16.8), we obtain
the nondegenerate
bilinear
pairing
$$
\roman H_{\roman{par}}^0(\pi,\{\pi_j\};\fra g_{\phi})
\otimes
\roman H_{\roman{par}}^2(\pi,\{\pi_j\};\fra g_{\phi})
@>>>
\Bobb R.
\tag2.16.9
$$
It has been shown
in Section 3 of \cite\guhujewe\ 
that $\cdot$
induces a symplectic structure
$$
\roman H_{\roman{par}}^1(\pi,\{\pi_j\};\fra g_{\phi})
\otimes
\roman H_{\roman{par}}^1(\pi,\{\pi_j\};\fra g_{\phi})
@>>>
\Bobb R
\tag2.16.10
$$
in degree 1.
\smallskip
These observations yield the infinitesimal structure
of the spaces 
$\roman{Hom}(F,G)_{\bold C}$ and
$\Cal H(\Cal P,G)_{\bold C}$ 
by means 
of the commutative diagram
$$
\CD
@.
\roman T_eG
@>\roman T_e\alpha_{\phi}>>
\roman T_{\phi} \left (\roman{Hom}(F,G)_{\bold C}\right)
@>{\roman T_{\phi} \rh}>>
\roman T_{r (\phi)}G
\\
@.
@A{\roman{Id}}AA
@A{\roman R_{\phi}}AA
@A{\roman R_{r(\phi)}}AA
\\
\bold C_{\roman{par}}(\Cal P, \fra g_{\phi})\colon
@.
\fra g
@>>{\delta^0_{\phi}}>
\fra g^{2\ell} \times \fra h_1 \times \dots\times \fra h_n
@>>{\delta^1_{\phi}}>
\fra g,
\endCD
\tag2.16.11
$$
having its vertical arrows isomorphisms
of vector spaces.  
This diagram
arises from the diagram (2.16.2)
by restriction.
\smallskip\noindent
(2.17) {\smc Proofs of Theorems {\rm 2.14} and {\rm 2.13}}.
The hypothesis
that the given
invariant symmetric bilinear form $\cdot$ on $\fra g$
be nondegenerate is still in force.
We now spell out explicitly an argument
already used in Section 5 of \cite\modus,
in Section 2 of \cite\modustwo, and in Section 8 of \cite\guhujewe.
\smallskip
For $1 \leq k \leq n$,
let $\tau_k$ be the 2-form
on $C_k$
constructed in \cite{\guhujewe~(6.2)}
and satisfying
$$
\roman{exp}^* \tau_k = \beta - \omega_k
\tag2.17.1
$$
on $\Cal O_k$.
We recall that
the equivariant  2-form
$\omega_{c,\wide {\Cal P},\bold C}$
on
$\Cal H(\wide {\Cal P},G)_{\bold C}$
given in \cite{\guhujewe~(7.1.1)}
is defined by
$$
\omega_{c,\wide {\Cal P},\bold C}
=
\eta^*\omega_c
-\widehat r^* \beta
+
\overline z_1^*\tau_1
+
\dots
+
\overline z_n^*\tau_n.
\tag2.17.2
$$

\proclaim {Lemma 2.17.3}
Near the zero locus
of the momentum mapping
$
\mu_{\bold C}
\colon
\Cal H({\Cal P},G)_{\bold C}
@>>>
\fra g_0^*
$
brought  into play in {\rm (2.14.6)},
that is, near the space $\roman{Hom}(\pi, G)_{\bold C}$, viewed as
a subspace of
$\Cal H(\Cal P, G)_{\bold C}$
or, more generally,
near the preimage 
with respect to $\mu_{\bold C}$
of the center of $\fra g$,
the 2-form $\omega_{c,\wide {\Cal P},\bold C}$
on $\Cal H(\Cal P, G)_{\bold C}$
is nondegenerate.
\endproclaim

\demo{Proof}
Let $K = G$ and $M=\Cal H(\Cal P, G)_{\bold C}$,
let $p$ be a point 
of $\Cal H(\Cal P, G)_{\bold C}$
having the property that the group (!) homomorphism 
$\phi =\eta (p)$ from $F$ to $G$,
cf. (2.14.4) for the notation $\eta$,
sends $r$ to a central element of $G$,
and 
apply Lemma 2.15.4,
with $\omega_{c,\wide {\Cal P},\bold C}$
playing the role of $\omega$.
To see that we are in the situation of that Lemma,
we proceed as follows.
Inspection of (2.16.11) shows that
right translation identifies
the vector space $H_p$ in (2.15.4)
with the parabolic cohomology group
$\roman H_{\roman{par}}^1(\pi,\{\pi_j\};\fra g_{\phi})$
and the form $\omega_p$
with the induced 2-form 
(2.16.10).
Moreover, the stabilizer Lie algebra $\fra k_p$
and the cokernel of
$d\mu_p\colon \roman T_p M \to \fra g^*$
amount to
the zero'th and second cohomology group, respectively,
of the chain complex
(2.16.4) defining parabolic cohomology
whence
$$
\fra k_p
\cong
\roman H_{\roman{par}}^0(\pi,\{\pi_j\};\fra g_{\phi}),
\quad
\roman{coker}(d\mu_p)
\cong 
\roman H_{\roman{par}}^2(\pi,\{\pi_j\};\fra g_{\phi}).
$$
The 
nondegeneracy of the resulting pairing
(2.16.9)
implies that
the image $d\mu_p(\roman T_pM)$ in $\fra g^*$
consists precisely of the linear forms on $\fra g$
which annihilate the stabilizer Lie algebra $\fra k_p$.
Furthermore,
when $p$ does not lie in the zero locus of the
momentum mapping,
adding a suitable constant, we may change the momentum
mapping so that $p$ lies in the zero locus
without changing the remaining geometrical data.
Hence the hypotheses of Lemma 2.15.4 are satisfied,
whence the form $\omega_{c,\wide {\Cal P},\bold C}$
is nondegenerate 
near the point $p$ as asserted. \qed \enddemo

\demo{Proof of {\rm 2.14}}
The mapping 
$i^{\sharp}$, cf. (2.14.5),
fits into the pull back square
$$
\CD
\Cal H(\widetilde {\Cal P},G)_{\overline{\Cal O}}
@>{(\widetilde z_1,\dots, \widetilde z_n)}>>
\Cal O_1 \times \dots \times \Cal O_n
\\
@V{i^{\sharp}}VV
@VV{\roman{exp} \times \dots\times \roman{exp}}V
\\
\Cal H({\Cal P},G)_{\bold C}
@>{(z_1,\dots,z_n)}>>
C_1 \times \dots \times C_n
\endCD
\tag2.14.7
$$
whose  vertical maps are the corresponding induced maps;
since the exponential maps,
restricted to the
$\Cal O_k$'s, are covering projections,
so is
the map $i^{\sharp}$.
Plainly,
then,
$$
(i^{\sharp})^*
\omega_{c, {\Cal P},\bold C}
=
(i^{\sharp})^*\eta^*\omega_c
-
(i^{\sharp})^*\widehat r^* \beta
+
\widetilde z_1^*(\beta - \omega_1)
+
\dots
+
\widetilde z_n^*(\beta - \omega_n).
\tag2.14.8
$$
Comparison with (2.11)
shows that
this is the 2-form
$\omega_{\widetilde c,\widetilde {\Cal P},\omega}$
on
$\Cal H(\widetilde {\Cal P},G)_{\overline{\Cal O}}$
to which the closed 2-form
$\omega_{\widetilde c,\widetilde {\Cal P}}$
on $\Cal H(\widetilde {\Cal P},G)$
descends by reduction from
$\Cal H(\widetilde {\Cal P},G)$
to
$\Cal H(\widetilde {\Cal P},G)_{\overline{\Cal O}}$.
More precisely, the constituents, respectively,
$$
{\widetilde \eta}^* \omega_{\widetilde c},\,
\widehat {\widetilde r}^* \beta,\,
\widehat a_1^* \beta,
\dots,
\widehat a_n^* \beta
$$
and
$$
(i^{\sharp})^*\eta^*\omega_c,\,
(i^{\sharp})^*\widehat r^* \beta,\,
\widetilde z_1^*\beta, 
\dots,
\widetilde z_n^*\beta  
$$
correspond to each other,
and the additional terms
$-\widetilde z_1^*\omega_1,\dots,-\widetilde z_n^*\omega_n$
arise since  reduction is applied to
$\Cal H(\widetilde {\Cal P},G)$
at the (co)adjoint orbit
$\Cal O_1 \times \dots \times \Cal O_n$
for
$\Gamma = G_1 \times \dots \times G_n$.
The distinction between the groupoid chain
$\widetilde c$ and the group chain
$c$ disappears
when the action of
$\Gamma$ is divided out,
since then the groupoid generators
$\gamma_k$ are no longer relevant.
By Lemma 2.17.3,
the 2-form $\omega_{c,\wide {\Cal P},\bold C}$
on
$\Cal H(\Cal P, G)_{\bold C}$
is nondegenerate 
near the zero locus.
Consequently
the 2-form
$\omega_{\widetilde c,\widetilde {\Cal P},\omega}$
on
$\Cal H(\widetilde {\Cal P},G)_{\overline{\Cal O}}$
is nondegenerate.
We have already observed
that the momentum mappings
(2.14.2) and (2.14.6)
correspond as well.
This proves Theorem 2.14. \qed
\enddemo

\demo{Proof of {\rm 2.13}}
Let $q \in \Cal H(\widetilde {\Cal P},G)$.
The image of $q$ under the momentum mapping
$$
\mu =\mu_1\times \dots \times \mu_n\colon
\Cal H(\widetilde {\Cal P},G)
@>>>
\fra g^*_1\times \dots \times \fra g^*_n
$$
lies in precisely one
(co)adjoint orbit
$\Cal O_1 \times \dots \times \Cal O_n$
for $\Gamma = G_1 \times \dots \times G_n$
and, by construction, 
each $O_k$ consists of regular points for the exponential map.
The reduced space
$\Cal H(\widetilde {\Cal P},G)_{\overline {\Cal O}}$
arises from reduction
at zero,
applied to 
the diagonal action of 
$\Gamma =G_1\times \dots \times G_n$
on 
$$
\Cal H(\widetilde {\Cal P},G)
\times
\overline{\Cal O}_1 \times \dots \times \overline{\Cal O}_n
$$
which is
hamiltonian, with momentum mapping
$$
\mu - \iota \colon
\Cal H(\widetilde {\Cal P},G)
\times
\overline{\Cal O}_1 \times \dots \times \overline{\Cal O}_n
@>>>
\fra g^*_1\times\dots\times \fra g^*_n,
$$
where
$\iota$
denotes
the inclusion of 
$\Cal O_1 \times \dots \times \Cal O_n$
into
$\fra g^*_1\times\dots\times \fra g^*_n$.
The resulting presymplectic structure
on
$\Cal H(\widetilde {\Cal P},G)
\times
\overline{\Cal O}_1 \times \dots \times \overline{\Cal O}_n$
is given by
$$
(\omega_{\widetilde c,\widetilde {\Cal P}},
-\omega_1,\dots,-\omega_n);
\tag2.13.1
$$
here 
$\omega_{\widetilde c,\widetilde {\Cal P}}$
is the 2-form (2.11.1).
Let
$$
p = (q,\mu(q))
\in M=
\Cal H(\widetilde {\Cal P},G)
\times
\overline{\Cal O}_1 \times \dots \times \overline{\Cal O}_n,
$$
and write
$[p] \in \Cal H(\widetilde {\Cal P},G)_{\overline {\Cal O}}$
for the point represented by $p$.
Then,
with reference to the notation introduced in (2.15),
$\omega_{c,\wide {\Cal P},\bold C}$
now playing the role of $\omega$,
the canonical map from
$H_p$ to
$\roman T_{[p]}\Cal H(\widetilde {\Cal P},G)_{\overline {\Cal O}}$
is an isomorphism
identifying the form $\omega_p$
with the reduced form
$\omega_{\widetilde c,\widetilde {\Cal P},\omega}$
(cf. (2.14.1))
at $[p]$.
A dimension count shows that, in this case,
the stabilizer Lie algebra $\fra k_p$ is trivial
and that the point $p$ is regular for the momentum mapping
$\mu - \iota$.
By (2.14),
the 2-form
$\omega_{\widetilde c,\widetilde {\Cal P},\omega}$
is nondegenerate at $[p]$.
Hence the hypotheses of (2.15.4)
are satisfied whence
the 
2-form
(2.13.1)
on
$\Cal H(\widetilde {\Cal P},G)
\times
\overline{\Cal O}_1 \times \dots \times \overline{\Cal O}_n$
is nondegenerate near
$p$ and hence
the 2-form
$\omega_{\widetilde c,\widetilde {\Cal P}}$
on $\Cal H(\widetilde {\Cal P},G)$
is nondegenerate near $q$.
This completes the proof of (2.13).\qed
\enddemo

\smallskip
\noindent
(2.18) {\smc The proof of main result}.
We now prove the theorem spelled out in the introduction.
\smallskip
Consider the smooth manifold
$\Cal H(\widetilde {\Cal P},G)$,
with 2-form
$\omega_{\widetilde c,\widetilde {\Cal P}}$,
cf. (2.11.1);
by (2.13),
near the zero locus of 
$\widehat{\widetilde r}\colon \Cal H(\widetilde {\Cal P},G) \to O_0$
or, more generally,
near the preimage 
of the center of $\fra g$,
this 2-form
is nondegenerate, that is, 
a symplectic structure.
Let $\Cal M \subseteq\Cal H(\widetilde {\Cal P},G)$
be an open $(G_0 \times \dots \times G_n)$-invariant submanifold
where
$\omega_{\widetilde c,\widetilde {\Cal P}}$
is nondegenerate
and which
contains the preimage 
of the center of $\fra g$.
Abusing the notation
introduced in (2.11), we shall
denote the restrictions of 
$\omega_{\widetilde c,\widetilde {\Cal P}}$,
$\mu$,
$\mu_0$,
$\mu_1$, etc.
to $\Cal M$  
by,
respectively,
$\omega_{\widetilde c,\widetilde {\Cal P}}$,
$\mu$,
$\mu_0$,
$\mu_1$, etc. as well.
The induced action of $G_0 \times \dots \times G_n$
on $\Cal M$ is hamiltonian, with momentum mapping
$$
\mu=\mu_0 \times \mu_1 \times \dots\times \mu_n
\colon
\Cal M 
@>>>
\fra g_0^* \times \fra g_1^*\times \dots \times \fra g_n^* .
$$
The data $(\Cal M,\omega_{\widetilde c,\widetilde {\Cal P}}, \mu)$
constitute the {\it extended moduli space\/}
which we are looking for.
\smallskip
By the main result of \cite\sjamlerm,
the reduced space 
$\Cal M_0 = \mu_0^{-1}(0)\big/ G_0$,
for the action
of the zero'th copy  $G_0$ of $G$ on $\Cal M$,
inherits a stratified symplectic space structure.
The action
of $G_0 \times G_1\times \dots \times G_n$ on $\Cal M$
manifestly induces a hamiltonian action 
(in the sense of stratified symplectic spaces)
of
$\Gamma =G_1\times \dots \times G_n$
on $\Cal M_0$.
The 
subalgebra 
$(C^{\infty}(\Cal M_0))^\Gamma$
of $\Gamma$-invariant
functions in 
$C^{\infty}(\Cal M_0)$
then yields a Poisson algebra 
$\left(C^{\infty}(\Cal M_0\big/ \Gamma),\{\cdot,\cdot\}\right)$
of continuous functions
on 
the quotient space $\Cal M_0\big/ \Gamma$.
\smallskip
Recall that the smooth map
$\widetilde \eta \colon \Cal H(\widetilde {\Cal P},G)
\to
\roman{Hom}(\widetilde F,G)$
has been brought into play in (2.4), and write
$\Cal H(\widetilde \pi,G) = \widetilde \eta(\mu_0^{-1}(0))$.
This is the $\widetilde \eta$-image
in 
$\roman{Hom}(\widetilde F,G)$
of 
the zero locus $\mu_0^{-1}(0)$ of the map $\mu_0$
from 
$\Cal H(\widetilde {\Cal P},G)$
to $O_0$;
the space $\Cal H(\widetilde \pi,G)$
clearly lies in
$\Cal M$
as the zero locus $\mu_0^{-1}(0)$ of $\mu_0$, now viewed 
as a map from 
$\Cal M$ 
to $O_0$.
Restriction
from the groupoid $\widetilde \pi$
to the group $\pi$, cf. (2.3),
induces a homeomorphism
$$
i^* \colon \roman{Rep}(\widetilde \pi,G)
=
\roman{Hom}(\widetilde \pi,G)/G^{n+1}
@>>>
\roman{Rep} (\pi,G)
$$
which passes to a homeomorphism
$$
i^\natural \colon 
\Cal H (\widetilde \pi,G)/G^{n+1}
@>>>
\Cal R (\pi,G).
$$
Hence
$\widetilde \eta$ induces a projection map from
$\Cal M_0 = \mu_0^{-1}(0)\big/ G$
onto
$\Cal H (\widetilde \pi,G)/G^{n+1}$
and hence a projection map
$$
q
\colon
\Cal M_0\big/ \Gamma
@>>>
\Cal R (\pi,G).
$$
\smallskip
We now suppose that the $O_j$'s, $1 \leq j \leq n$, are all  equal to
the chosen
invariant connected neighborhood of zero, $O$,
where the exponential map
from $\fra g$ to $G$
is a diffeomorphism. Then
the map $\widetilde \eta$
is injective and $q$ is a homeomorphism.
Thus the Poisson algebra
$\left(C^{\infty}(\Cal M_0\big/ \Gamma),\{\cdot,\cdot\}\right)$
then yields a Poisson algebra 
$\left(C^{\infty}(\Cal R (\pi,G),\{\cdot,\cdot\}\right)$
of continuous functions
on
$\Cal R (\pi,G)$;
this is the Poisson algebra
which we are looking for.
Comparison of the construction of
this Poisson structure
with that of
$\left(C^{\infty}(\roman{Rep} (\pi,G)_{\bold C},\{\cdot,\cdot\}\right)$
in \cite\guhujewe\ 
shows that,
on each 
moduli space $\roman{Rep}(\pi,G)_{\bold C}$ (that
lies in $\Cal R(\pi,G)$),
the Poisson structure on
$\Cal R (\pi,G)$
restricts to the 
stratified symplectic 
Poisson algebra
on $\roman{Rep}(\pi,G)_{\bold C}$.
The proof of the theorem spelled out in the introduction
is now complete
save that
the additional statement saying that, for $G$ compact,
the subspace
$\Cal R(\pi,G)$ may be taken
dense in $\roman{Rep}(\pi,G)$,
has not been justified yet.
This will be achieved in the next section.

\smallskip
\noindent
{\smc Remark 2.19.}
The chosen neighborhood $O$ of the
origin of the Lie algebra $\fra  g$,
and hence its image
$B$ in $G$,
inherits a Poisson
structure from $\fra g$
(or rather from $\fra g^*$),
and this structure
cannot be extended to $G$.
It would be interesting to know
whether the obstruction
to extending
the Poisson structure on
$\Cal R (\pi,G)$
to one on
$\roman{Rep} (\pi,G)$
is formally of the same kind.
The simplest case occurs perhaps for
$\ell = 0$ and $n=3$.
(When
$\ell = 0$ and $n=1$ or $n=2$
the space
$\roman{Rep} (\pi,G)_{\bold C}$
is either empty
or consists of a single point.)
Then 
$\roman{Rep} (\pi,G)_{\bold C}$
is the space of $G$-orbits of triples
$(a,b,c) \in C_1 \times C_2 \times C_3$
satisfying $abc = e$ while
$\roman{Rep} (\pi,G)$
is the space of $G$-orbits of triples
$(a,b,c) \in G \times G \times G$
satisfying $abc = e$.
Even though
the induced Poisson structure on 
$B \subseteq G$
cannot be extended to a Poisson structure on $G$
it may still be possible to extend
the Poisson structure
on $\Cal R (\pi,G)$
to one on
$\roman{Rep} (\pi,G)$.
Another special case worth looking at arises with $G= \roman{SU}(2)$
and a surface $\Sigma$ 
of genus $\ell>0$ (say)
with a single boundary circle.
In this case,
when $O$ is taken to be the open ball in
$\fra g = \roman{su}(2)$ consisting of all $X$ such that
$-\psi(X,X) < 8 \pi^2$ where
$\psi$ is the Killing form,
the image $B$ of $O$ under the exponential map
is the open ball
$\roman{SU}(2) \setminus \{-\roman{Id}\} $
in $\roman{SU}(2)$, and
$\roman{Rep} (\pi,G)$ is the disjoint union of
$\Cal R (\pi,G)$ and $\roman{Rep}_{-1} (\overline \pi,G)$;
here
$\roman{Rep}_{-1} (\overline \pi,G)$
refers to the space of {\it twisted\/} representations
$\phi$ of the fundamental group
$\overline \pi$
of the corresponding closed surface 
$\overline \Sigma$
having the property that,
after a choice $x_1,y_1, \dots, x_\ell, y_\ell$
of generators for $\overline \pi$ has been made,
$$
[\phi x_1,\phi y_1] \cdot \dots \cdot [\phi x_\ell,\phi y_\ell] =
-\roman{Id} \in G= \roman{SU}(2).
$$
The construction in the present paper
endows
$\Cal R (\pi,G)$
with a Poisson algebra
of continuous functions, and
the space
$\roman{Rep}_{-1} (\overline \pi,G)$
has been known 
for a while
to carry a structure of a smooth symplectic manifold
and hence in particular is endowed with the corresponding
symplectic Poisson algebra structure
\cite{\atibottw, \guhujewe, \modus, \huebjeff, \jeffrtwo,
\narasesh}.
In the vector bundle picture, 
after a  choice of (compatible) complex structure on
$\overline \Sigma$
has been made,
$\roman{Rep}_{-1} (\overline \pi,G)$ 
may be identified with
the space
of stable holomorphic
 vector bundles on 
$\overline \Sigma$
having rank 2, degree 1, and trivial determinant.
(There is no difference between semistable and stable vector bundles
in this case.)
But is there a Poisson algebra
of continuous functions
on 
$\roman{Rep} (\pi,G)$ which, on
$\Cal R (\pi,G)$ and $\roman{Rep}_{-1} (\overline \pi,G)$,
restricts to the corresponding Poisson algebras?
\smallskip
For general $\ell$ and $n$ and a general group $G$,
the difference
$\roman{Rep} (\pi,G) \setminus \Cal R (\pi,G)$
is the union 
in $\roman{Rep} (\pi,G)$
of all 
spaces $\roman{Rep}(\pi,G)_{\bold C}$
for those 
$n$-tuples 
$\bold C = (C_1,\dots,C_n)$
of conjugacy classes which are {\it not $B$-regular\/} in the sense
that at least one of the $C_k$ is not $B$-regular.
However,
in view of the main result of \cite\guhujewe, 
the spaces 
$\roman{Rep}(\pi,G)_{\bold C}$
also arise by symplectic reduction, applied  to suitable extended moduli spaces
and hence inherit structures of stratified symplectic spaces.
Among the  obstacles
to extending the Poisson structure
from
$\Cal R (\pi,G)$ to $\roman{Rep} (\pi,G)$
is the following one:
The space
$\Cal H(\widetilde {\Cal P},G)$
coming into play in
(2.4) is contained in a larger space
$\roman H(\widetilde {\Cal P},G)$
to be defined by means of a pull back square
$$
\CD
\roman H(\widetilde {\Cal P},G)
@>{(\widehat {\widetilde r},\widehat a_1,\dots,\widehat z_n)}>>
\fra g_0 \times \fra g_1\times\dots\times \fra g_n
\\
@V{{\widetilde \eta}}VV
@VV{\roman{exp}  \times \dots\times \roman{exp}}V
\\
\roman{Hom}(\widetilde F,G)
@>>{(\widetilde r,a_1,\dots,a_n)}> G_0\times G_1 \times \dots \times G_n.
\endCD
$$
However,
$\roman H(\widetilde {\Cal P},G)$
cannot naively serve as an extended moduli space
since
it is {\it not\/} a {\it smooth\/} manifold:
its singularities arise from the non-regular points of the
exponential map.
Moreover the relationship
between the zero locus of
$\widehat {\widetilde r}$
and
$\roman{Rep} (\pi,G)$
is rather intricate and does not just come down
to 
a homeomorphism
after  the action of
$\Gamma$ on the zero locus
has been divided out.

\smallskip
\noindent
{\smc Remark 2.20.}
In \cite{\jeffrone \ (5.2)},
a certain space 
denoted $\Cal N^{\fra g,N}$
has been constructed
and, by gauge theory methods,
a symplectic structure
on a smooth open part of
$\Cal N^{\fra g,N}$
together with a hamiltonian action
of a product $G^N$ of finitely many copies of $G$
on this space and a momentum mapping into
$(\fra g^*)^N$
have been obtained.
An appropriate smooth open part of the space
$\Cal N^{\fra g,N}$, 
where $N=n$, arises from
$\Cal H(\widetilde {\Cal P},G)$
by reduction with respect to the momentum mapping
$\mu_0$.
We do not know whether
the symplectic structure
on $\Cal N^{\fra g,N}$ constructed
in \cite{\jeffrone \ (5.2)}
coincides with the symplectic structure
obtained above 
by a purely
finite dimensional construction 
via reduction from (2.11.1).
An investigation
of the difference between adjoint orbits and conjucacy classes,
crucial for the construction of the Poisson structure
on the ambient space,
has not been undertaken 
in \cite{\jeffrone}, though.

\beginsection 3. The case when $G$ is compact and connected

Throughout this section,
$G$ will be a compact and connected Lie group.
Our aim is to show that, roughly speaking,
there is a \lq\lq nice\rq\rq\ 
choice for the open neighborhood $O$ of the origin in
the Lie algebra $\fra g$ of $G$
having dense image $B$ in $G$ under the exponential map.
\smallskip
Recall that a point of $G$ is called {\it singular\/}
if it lies in two distinct maximal tori
and {\it regular\/} otherwise,
cf. e.~g.
\cite{\broediec~p.~168}.
The set of regular points is denoted by $G_r$.
The set $G \setminus G_r$
of singular points has codimension $\geq 3$
\cite{\broediec~(V.2.6)} and hence $G_r$ is, in particular, dense
in $G$.
Thus, in view of the well known structure
of $G_r$
for $G$ compact, connected, and simply connected,
(and hence semisimple) \cite{\broediec~(V.7.11)}, 
we could in fact choose
$O$ in $\fra g$ in such a way that
the image $B=\roman{exp}(O)$ coincides with $G_r$.
With this choice, we would immediately get
$\Cal R (\pi,G)$ dense in
$\roman{Rep}(\pi,G)$.
This choice of $O$ (and hence $B$) will not be \lq\lq nice\rq\rq, though,
since $B$ would not then contain e.~g. the neutral element of $G$.
We now show how $G_r$ may be enlarged to a suitable subset $B$ of $G$
having better properties for our purposes.
\smallskip
Choose a root of the equation $z^2+1=0$ and denote it by $i$ as usual.
Let  $O$
be the neighborhood of the origin in $\fra g$
consisting of those $X$ in $\fra g$
which have the property that the endomorphism $\roman{ad}(X)$
of $\fra g$ has only eigenvalues $\lambda = 2 \pi i \nu$
with $|\nu| <1$.
The following is presumably known;
we did not find it in the literature.

\proclaim{Lemma 3.1}
When $G$ is and simply connected (and hence semisimple),
the 
exponential map, restricted to $O$,
is a diffeomorphism
onto its image $B$ in $G$.
\endproclaim

We now prepare for the proof of this Lemma.
To fix notation, we recall that, after a choice
of maximal torus $T$ in $G$ has been made,
the global weights $\vartheta \colon T \to S^1$
of the adjoint representation are called
{\it global roots\/},
the
corresponding linear forms $\Theta \colon T \to i \Bobb R$
the {\it infinitesimal roots\/},
and the real linear forms
$\alpha \colon T \to \Bobb R$
such that
$\Theta = 2 \pi i \alpha$
are said to be the {\it real roots \/} of $G$.
Given a point $X \in \fra g$ which is singular for
the exponential map,
$\roman{exp}(X)$ is a singular
point of $G$
but the converse is not necessarily true:

\proclaim{Proposition 3.2}
A point $X$ of $\fra g$ is singular for
the exponential map
if and only if 
$\roman{exp}(X)$ is a singular point of $G$
(in the sense of Lie groups)
in such a way that,
for some real root $\alpha$ and some integer $k \ne 0$,
$\alpha (X) = k$.
\endproclaim

\demo{Proof}
Let $X \in \fra g$;
it lies in a Cartan subalgebra of $\fra g$, and
the point $\roman{exp}(X)$ lies in the corresponding maximal torus $T$.
With respect to a basis of weight vectors,
the complexification of $\roman{ad}(X)$ is then in diagonal form,
with entries
$\Psi(X)$ and 0
where
$\Psi$ runs through the infinitesimal roots.
Now
$\roman{exp}(X)$ is a singular point of $G$
if and only if,
for some real root $\alpha$
and some integer $k$ (which may be zero),
$\alpha (X) = k$.
On the other hand, the exponential map is
non-regular at a point $X$ of the Lie algebra $\fra g$
if and only if the endomorphism
$\roman{ad}(X)$ of $\fra g$ has an eigenvalue of the kind
$2 \pi i k, k \ne 0$,
$k$ being a nonzero integer,
cf. Corollary 1 on p. 106 of \cite{\kirilboo}.
This implies the assertion. \qed
\enddemo

We now spell out two immediate consequences of (3.2).

\proclaim {Corollary 3.3}
The subset $O$ is the
(uniquely determined)
largest {\rm connected\/} neighborhood $O$ of the origin
of $\fra g$ where the exponential map is regular. \qed
\endproclaim

Write $O_r$ for the subset of $O$
consisting of those $X$ in $\fra g$
which have the property that $\roman{ad}(X)$
has only eigenvalues $\lambda = 2 \pi i \nu$
with $0 < |\nu| <1$.

\proclaim{Corollary 3.4}
When $G$ is simply connected,
the subset $O_r$ is mapped diffeomorphically onto the set $G_r$
of regular points of $G$. \qed
\endproclaim

\smallskip
We choose a maximal torus
$T$ in $G$, with Lie algebra $\fra t$;
let $r = \dim T (= \roman{rank} (G)$).
The {\it diagram\/}
is the inverse image in $\fra t$ of the set
of singular points (in the sense of Lie groups) in $T$,
with respect to the exponential map
\cite{\broediec~(V.7)}.
Alternatively, in view of (3.2),
the diagram  consists of the hyperplanes
$L_{\alpha n} = \alpha^{-1}(n)$ in $\fra t$
where $\alpha$ runs through the real roots of $G$ and
$n$ through the integers.
Further, a choice of {\it fundamental Weyl chamber\/} corresponds 
to a choice of a 
system $S$ of {\it simple\/} 
roots and vice versa; the {\it walls\/} of the Weyl chamber are 
then the hyperplanes 
$L_{\alpha} =L_{\alpha 0}$
where $\alpha$ runs through the corresponding simple roots.
\smallskip
We now choose a fundamental Weyl chamber; 
in terms of the corresponding system $S$ of simple (real) roots,
this Weyl chamber consists of those
$X$ in $\fra t$ which satisfy the condition
$\alpha (X) >0$ for every $\alpha \in S$.
The choice of fundamental Weyl chamber uniquely determines
an  alcove  (in $\fra t$), that is, a connected component 
$P$
of the complement
of the diagram,
in such a way that
the origin of $\fra t$ is a vertex of the closure
$\overline P$
and that,
for each wall of the fundamental Weyl chamber, a certain 
(uniquely determined)
convex subset thereof
containing the origin of $\fra t$
constitutes a wall of
$P$ (but $P$ will have other walls).
The map from $(G/T) \times P$ to $G$ which assigns
$x \roman {exp}(\Lambda) x^{-1} $ to
$(xT,\Lambda) \in (G/T) \times P$
is known to be a universal covering
with a connected
total space onto the subset $G_r$ of regular elements
of $G$,
in fact, a diffeomorphism
when $G$ is simply connected,
see e.~g. \cite{\broediec~(V.7.11)}.
\smallskip
Let
$\widetilde P$ be the subspace
of the closed alcove
$\overline P$ 
which is obtained when those walls of
$\overline P$ 
are removed which
do not
lie in any of the walls of the chosen Weyl chamber.
In other words,
$$
\widetilde P
=\{X \in \fra t; \ \alpha (X) \geq 0\  \text{for} \  \alpha \in S,
\ 
\alpha (X) < 1\  \text{for} \  \alpha \not\in S\}.
$$
Since, for any $X \in \fra g$,
with respect to a basis of weight vectors, the complexification
of $\roman{ad}X$ is in diagonal form, with entries
$\pm 2 \pi i \alpha(X)$ and $0$ where $\alpha$ runs through the positive roots
of $G$, the subspace
$\widetilde P$ of
$\overline P$
consists precisely of those $X$ in $\fra t$
which have the property that the endomorphism $\roman{ad}(X)$
of $\fra g$ has only eigenvalues $\lambda = 2 \pi i \nu$
with $|\nu| <1$.
Consequently 
the image of the map from $(G/T) \times \widetilde P$ to $\fra g$ which assigns
$\roman{Ad}(x) \Lambda $ to
$(xT,\Lambda) \in (G/T) \times \widetilde P$
is precisely
the neighborhood $O$ of the origin
in $\fra g$.
In view of
(3.2),
every point of
this neighborhood $O$ 
is regular for the exponential map. 
\smallskip

\demo{Proof of Lemma {\rm 3.1}}
Suppose that $G$ is simply connected.
We claim that then
the restriction of the exponential map to $O$
is injective.
It is clearly injective when restricted to
the subset $O_r$
since
the map 
from
$G/T \times P$ to $G$
which sends
$(xT,\Lambda) \in (G/T) \times P$
to
$x \roman {exp}(\Lambda) x^{-1} $ 
is a diffeomorphism onto
the subset $G_r$ of regular elements
of $G$.
Furthermore, $G$ being
simply connected,
the map from $\overline P$ to $G$ which maps
$\Lambda$ to
$\roman {exp}(\Lambda)$
yields a bijective correspondence
between the points of
$\overline P$
and
the conjugacy classes of $G$ and hence
between adjoint orbits
and conjugacy classes of $G$.
We now claim that
this implies the injectivity 
of the restriction of the exponential map to $O$.
In fact,
let $X,Y$ be two points of $O$
having the same image under the exponential map.
If one of them is in $O_r$,
the other one must be in $O_r$, too,
since $O_r$ is the preimage of the regular points of $G$,
and hence $X$ and $Y$ must coincide, 
since the corresponding map from $(G/T) \times O_r$
to $G_r$ is a diffeomorphism.
If the points $X,Y$ of $O$ are both not
in $O_r$,
they necessarily lie in the same adjoint orbit
$\Cal O \subseteq O$
since
the map from $\overline P$ to $G$ 
sending $\Lambda \in \overline P$ to
$\roman{exp}\Lambda$
yields a bijective correspondence
between adjoint orbits
and conjugacy classes of $G$.
Now, the exponential map, still being regular
on $\Cal O$,
is just a covering projection
from 
$\Cal O$
onto its image, which is a conjugacy class (say) $C$.
However, 
$G$ being simply connected,
in view of an unpublished result
of Bott's,
cf. \cite\borelone\ (Theorem 3.4)
and \cite\rasheone, 
the conjugacy class $C$ is simply connected or, equivalently,
the centralizer $Z_x$ of any
$x \in G$ is connected;
see also
what is said on p. 351 of \cite\helgaboo.
Hence
the covering projection from
$\Cal O$ onto $C$ is a diffeomorphism.
Consequently
the exponential map, restricted to
$\Cal O$, is injective. \qed
\enddemo

Next we consider a compact, connected, semisimple Lie group $G$
which is not necessarily simply connected.
Then the fundamental group $\pi_1G$ acts on the ball
$\widetilde B$ of regular values of the exponential map
from the Lie algebra $\fra g$ to the universal cover $\widetilde G$,
and we may choose for $B$ the image in $G$ under the covering projection
of a suitable fundamental domain
for the $\pi_1G$-action on $\widetilde B$.
We now describe this somewhat more explicitly:
The closed alcove $\overline P$
is a regular polyhedron, (referred to sometimes as {\it fundamental
polyhedron\/} in the literature,) and
the fundamental group $\pi_1G$
may be realized as a finite group of automorphisms thereof,
cf. \cite{\broediec~(V.7.17)~Ex.~5}.
With reference to this group of automorphisms,
we then choose a fundamental domain $D$ in $P$
whose closure $\overline D$
contains the origin and meets each wall of the Weyl chamber
in a convex subset of dimension $r-1$.
Let $\widetilde D$
be the subspace
of $\overline D$ 
which is obtained when those walls of
$\overline D$ 
are removed which
do not
lie in any of the walls of the chosen Weyl chamber,
and let $O$ be the image of
$\overline D$ in $\fra g$ under the canonical map
from $G/T \times \overline D$ to $\fra g$;
this set $O$ is a neighborhood of the origin in $\fra g$.
Further, the exponential map, restricted to $O$,
is a diffeomorphism from
$O$ onto its image $B$ in $G$, and $B$ is a neighborhood of the
neutral element of $G$ which contains $G_r$ by construction and
hence is, in particular, dense in $G$.
For $G=\roman{SO}(3)$,
the closed alcove
$\overline P$
is a line segment 
having the origin $o$ of $\fra t \cong \Bobb R$
as a vertex,
and 
$\overline D = \overline {oQ}$,
$Q$ being the midpoint of $P$,
is that half of
$\overline P$
which still contains
the origin $o$
as a vertex;
then $\widetilde D = 
\overline D \setminus \{Q\}$,
and the point $Q$ corresponds to the conjugacy class
of
$\left[\matrix -1& 0 & 0\\
                0&-1 & 0\\
                0& 0 & 1
       \endmatrix
 \right]$;
this explains
the remark made about 
$\roman{SO}(3)$
in the introduction.
\smallskip
Finally, let $G$ be an arbitrary compact connected Lie group.
Then a suitable covering $\widetilde G$ of $G$
may be written as a product
$\widetilde G = T \times K$
where $T$ is a torus and $K$ is compact and 
simply connected.
We may then 
apply to $K$ what has been said already  in order to find
an open neighborhood $O_K$ of the origin of the Lie algebra
$\fra k$ of $K$ and, furthermore, 
we may choose an open subset $O_T$
of the origin of the Lie algebra
$\fra t$ of $T$,
in such a way that
the exponential maps are diffeomorphisms
onto their images $B_T$ and $B_K$ in $T$ and $K$, respectively,
and that these images are dense in the corresponding groups.
A suitable choice of fundamental domain
in $O_T \times O_K$ for the action of the group
of deck transformations of the covering
from
$T \times K$
onto $G$
will then yield 
an open neighborhood $O$ of the origin of the Lie algebra
$\fra g$ of $G$ such that
the exponential map, restricted to $O$, 
is a diffeomorphism onto its image $B$ in $G$,
and so that $B$ is, furthermore, dense in $G$.
\smallskip
These considerations show that, for any compact and connected
Lie group $G$,
the subspace $\Cal R(\pi, G)$
may be taken dense in
$\roman{Rep}(\pi,G)$
in a certain \lq\lq nice way\rq\rq\ 
made precise above.
The proof of the theorem in the introduction is now complete.

\smallskip
\noindent
{\smc Example 3.5.} We give a brief description of
adjoint orbits and conjugacy classes
and their mutual relationships,
for the group $\roman {SU}(n)$;
since the other compact semisimple linear groups appear as closed
subgroups of the special unitary group this gives certain information
about other groups for free.
\smallskip
The group $\roman {SU}(n)$ consists of complex $(n \times n)$-matrices $A$ 
of determinant 1
such that
$\overline A^t A = E$, 
the standard maximal torus being the subgroup
of diagonal matrices $\roman{diag}(\zeta_1,\dots, \zeta_n)$
with $\zeta_1 \cdot \dots \cdot \zeta_n = 1$
where $\zeta_j = \roman{exp}(2 \pi i \nu_j)$, for $ 1 \leq j \leq n$.
Its Lie algebra
$\fra {su}(n)$
is the space of
skew-hermitian matrices of trace zero.
We write such a matrix as $M=2 \pi i N$;
the adjoint orbit through such an $M$ is determined by its eigenvalues 
$\lambda_1= 2 \pi i \nu_1,\dots,\lambda_n= 2 \pi i \nu_n$.
The Lie algebra $\fra t$ of the maxmial torus 
is the space of matrices
$M=2 \pi i N$
where $N = \roman{diag}(\nu_1, \dots, \nu_n)$
with
$\nu_1 + \dots + \nu_n= 0$.
The  simple real roots $\alpha_j$ are given by
$\alpha_j =\nu_j - \nu_{j+1}$, for $1 \leq j < n$,
and the fundamental Weyl chamber consists of
the matrices
$2 \pi i N$
where $N = \roman{diag}(\nu_1, \dots, \nu_n)$
with
$\nu_1 \geq \dots \geq \nu_n$.
The subspace $\widetilde P$ of the closed alcove
$\overline P$ 
introduced above
is then 
given by
the matrices
$2 \pi i \,\roman{diag}(\nu_1, \dots, \nu_n)$
in the fundamental Weyl chamber
which satisfy the additional condition
$\nu_1 - \nu_n<1$.
(The root $\nu_1 - \nu_n$ is the highest weight
of the adjoint representation.)
\smallskip
The stabilizer $Z_M$ in
$\roman{SU}(n)$
of $M=2 \pi i N$ in $\fra t$
depends on 
how many coincidencences there are amongst the
eigenvalues. 
Thus, if the first $n_1$ are equal,
the next $n_2$ are equal, and so on,
$$
Z_M = \left(\roman U(n_1) \times \roman U(n_2) \times \dots \times 
\roman U(n_s)\right)
\cap \roman{SU}(n).
$$
Then, if $M$ lies in the fundamental Weyl chamber,
it lies in
$n_1+ n_2+ \dots + n_s -s $ walls thereof, with
the convention that 
this number to be zero means that
$M$ lies in the interior 
(which only happens when  all the $n_j$ equal 1
and $s=n$).
The adjoint orbit $\Cal O_M$ through $M$ is then the
homogeneous space $\roman {SU}(n)/ Z_M$.
When $M$ lies in
$\widetilde P$,
the exponential map
identifies
$\Cal O_M$
with the conjugacy class $C_m$ through
$m=\roman{exp}(M) \in \roman {SU}(n)$.
However, when
$M$ lies in
$\overline P \setminus \widetilde P$,
the stabilizer $Z_m$ of 
$m=\roman{exp}(M) \in \roman {SU}(n)$
is strictly larger than
the stabilizer $Z_M$ of 
$M$
and the resulting fibre bundle projection
from 
the adjoint orbit $\Cal O_M$
onto 
the conjugacy class
$C_m$ through
$m$
has as fibre the homogeneous space
$Z_m\big/ Z_M$ which
has strictly  positive dimension.
\smallskip
For example, when $n=3$,
the fundamental Weyl chamber
is given by the triples $(\nu_1,\nu_2,\nu_3)$
with 
$\nu_1 \geq \nu_2 \geq\nu_3$ and
$\nu_1+\nu_2+\nu_3 = 0$;
its two walls consist of the points
$2\pi i\,\roman{diag}(\nu,\nu,-2\nu)$ and
$2\pi i\,\roman{diag}(2\nu,-\nu,-\nu)$, where $\nu \geq 0$.
The alcove $P$ is determined by the additional condition
$\nu_1 - \nu_3 <1$.
The conjugacy class through a regular point 
of $\roman {SU}(3)$
is the
six dimensional flag manifold $\roman {SU}(3)\big/(S^1 \times S^1)
(\cong \roman {U}(3) \big /(S^1)^3)$.
The points of
$\widetilde P$
which lie in any of the two walls of the Weyl chamber
look like
$M=2\pi i\,\roman{diag}(\nu, \nu, -2\nu)$
or
$M=2\pi i\,\roman{diag}(2\nu, -\nu, -\nu)$
with $0\leq \nu <\frac 13$;
for $\nu >0$,
the adjoint orbit $\Cal O_M$ 
through such a point $M$ 
is a 2-dimensional complex projective space
$P_2\Bobb C$,
and the exponential map induces a diffeomorphism onto the corresponding
conjugacy class.
However,
the points of
$\overline P \setminus \widetilde P$
are precisely of the kind
$M=2\pi i\,\roman{diag}(\nu, 1-2 \nu, \nu -1)$ 
with $\frac 13 \leq \nu \leq \frac 23$. 
In this case, for $\frac 13 < \nu < \frac 23$,
the conjugacy class $\Cal O_M$
through such a point 
$M$ is still a six dimensional flag manifold
of the kind
$\roman {U}(3) \big /(S^1)^3$
while
the conjugacy class $C_m$ through
$m = \roman{exp} M$
is a 2-dimensional complex projective space,
the projection mapping from
$\Cal O_M$
to
$C_m$
being the canonical one
from the flag manifold to
$P_2\Bobb C$
having as fibre 
a 1-dimensional complex projective space.
In the extreme cases
$\nu = \frac 13$ and
$\nu = \frac 23$, we get here the two points
$M_1=2\pi i\,\roman{diag}(\frac 13, \frac 13,-\frac 23)$ 
and
$M_2=2\pi i\,\roman{diag}(\frac 23, -\frac 13,-\frac 13)$,
respectively,
each of which lies in a wall of the fundamental Weyl chamber
as a vertex of $P$
but, beware, $M_1$ and $M_2$ do not lie in
$\widetilde P$;
the adjoint orbits $\Cal O_{M_1}$ and
$\Cal O_{M_2}$ of these points
are 2-dimensional complex projective spaces
but
their images $m_1$ and $m_2$, respectively, in $G$ under the exponential map
are central elements of $G$ and hence
the conjugacy classes of these images
consist each of 
a single point.

\medskip
\centerline{References}
\widestnumber\key{999}
\smallskip
\ref \no \atiyboo
\by M. F. Atiyah 
\book The geometry and physics of knots
\publ Cambridge University Press
\publaddr Cambridge, U. K.
\yr 1990
\endref
\ref \no  \atibottw
\by M. F. Atiyah and R. Bott
\paper The Yang-Mills equations over Riemann surfaces
\jour Phil. Trans. R. Soc. London  A
\vol 308
\yr 1982
\pages  523--615
\endref

\ref \no \bodenhu
\by H. U. Boden and Y. Hu
\paper Variations of moduli in parabolic bundles
\jour Math. Ann.
\vol 301
\yr 1995
\pages 539--559
\endref

\ref \no \borelone
\by A. Borel
\paper Sous-groupes commutatifs et torsion des groupes de 
Lie compacts connexes
\jour T\hataccent ohoku Math. J.
\vol 13
\yr 1961
\pages  216--240
\endref

\ref \no \broediec
\by T. Br\"ocker and T. tom Dieck
\book Representations of compact Lie groups
\bookinfo Graduate texts in Mathematics, Vol. 98
\publ Springer
\publaddr Berlin $\cdot$ Heidelberg $\cdot$ New York $\cdot$ Tokyo
\yr 1985
\endref

\ref \no \brownboo
\by R. Brown
\book Elements of modern topology
\publ McGraw-Hill
\publaddr London
\yr 1968
\endref
\ref \no \guhujewe
\by K. Guruprasad, J. Huebschmann, L. Jeffrey, and A. Weinstein
\paper Group systems, groupoids, and moduli spaces
of parabolic bundles
\jour Duke Math. J.
\vol 89
\yr 1997
\pages 377--412
\endref

\ref \no \helgaboo
\by S. Helgason
\book Differential Geometry, Lie Groups, and Symmetric Spaces
\bookinfo Pure and Applied Mathematics, 
A series of monographs and textbooks
\publ Academic Press
\publaddr New York $\cdot$ London $\cdot$ Toronto 
$\cdot$ Sydney $\cdot$ San Francisco
\yr 1978     
\endref

\ref \no \modus
\by J. Huebschmann
\paper Symplectic and Poisson structures of certain moduli spaces
\jour Duke Math. J.
\vol 80
\yr 1995
\pages 737--756
\endref
\ref \no \modustwo
\by J. Huebschmann
\paper Symplectic and Poisson structures of certain moduli spaces. II.
Projective
representations of cocompact planar discrete groups
\jour Duke Math. J.
\vol 80
\yr 1995
\pages 757--770
\endref
\ref \no \srni
\by J. Huebschmann
\paper
Poisson geometry of certain
moduli spaces
\paperinfo
Lectures delivered at the \lq\lq 14'th Winter School\rq\rq, Srni,
Czeque Republic,
January 1994
\jour Rendiconti del Circolo Matematico di Palermo, Serie II
\vol 39
\yr 1996
\pages 15--35
\endref
\ref \no \huebjeff
\by J. Huebschmann and L. Jeffrey
\paper Group cohomology construction of symplectic forms
on certain moduli spaces
\jour Int. Math. Research Notices
\vol 6
\yr 1994
\pages 245--249
\endref
\ref \no \jeffrone
\by L. Jeffrey
\paper 
Extended moduli spaces of flat connections
on Riemann surfaces
\jour Math. Ann.
\vol 298
\yr 1994
\pages 667--692
\endref
\ref \no \jeffrtwo
\by L. Jeffrey
\paper Symplectic forms on moduli spaces
of flat connections on 2-manifolds
\paperinfo in {\it Proceedings
of the Georgia International Topology Conference}, Athens, Ga.
1993, ed. by W. Kazez
\jour AMS/IP Studies in Advanced Mathematics
\vol 2
\yr 1997
\pages 268--281
\endref
\ref \no \karshone
\by Y. Karshon
\paper 
An algebraic proof for the symplectic
structure of moduli space
\jour Proc. Amer. Math. Soc.
\vol 116
\yr 1992
\pages 591--605
\endref
\ref \no \kirilboo
\by A. A. Kirillov
\book Elements of the theory of representations
\bookinfo Grundlehren der Mathematik, vol. 220
\publ Springer
\publaddr Berlin $\cdot$ Heidelberg $\cdot$ New York
\yr 1976
\endref

\ref \no \mehtsesh
\by V. Mehta and C. Seshadri
\paper Moduli of vector bundles on curves with parabolic structure
\jour Math. Ann.
\vol 248
\yr 1980
\pages 205--239
\endref

\ref \no \narasesh
\by M. S. Narasimhan and C. S. Seshadri
\paper Stable and unitary vector bundles on a compact Riemann surface
\jour Ann. of Math.
\vol 82
\yr 1965
\pages  540--567
\endref

\ref \no \rasheone
\by Rashevski
\paper A theorem on the connectedness of a subgroup of a simply 
connected Lie group commuting with any of its automorphisms
\jour Trans. Moscou Math. Soc.
\vol 30
\yr 1974
\pages 1--24
\endref

\ref \no \sjamlerm
\by R. Sjamaar and E. Lerman
\paper Stratified symplectic spaces and reduction
\jour Ann. of Math.
\vol 134
\yr 1991
\pages 375--422
\endref

\ref \no \trottone
\by H. F. Trotter
\paper Homology of group systems with applications to knot theory
\jour Annals of Math.
\vol 76
\yr 1962
\pages 464--498
\endref
\ref \no  \weiltwo
\by A. Weil
\paper  Remarks on the cohomology of groups
\jour                                      
Ann. of Math.
\vol 80                                    
\yr 1964
\pages  149--157
\endref

\ref \no  \weinsthi
\by A. Weinstein
\paper The symplectic structure on moduli space
\bookinfo
in: The Andreas Floer Memorial Volume,
H. Hofer, C. Taubes, A. Weinstein, and
E. Zehnder, eds;
Progress in Mathematics, Vol. 133
\publ Birkh\"auser
\publaddr Boston $\cdot$ Basel $\cdot$ Berlin
\yr 1995
\pages 627--635
\endref
\enddocument